\documentclass[12pt]{amsart}

\usepackage{enumitem}
\usepackage{psfrag}
\usepackage{graphicx}
\usepackage{pinlabel}
\usepackage{graphicx}
\usepackage{amsmath}
\usepackage{amssymb}
\usepackage{amscd}
\usepackage{autobreak}
\usepackage{picinpar}
\usepackage{times}
\usepackage{pb-diagram}
\usepackage{graphicx}
\usepackage{wrapfig}
\usepackage{xspace}
\usepackage{hyperref}
\usepackage{url}
\usepackage{caption}
\usepackage{subcaption}
\usepackage{fancyhdr}
\usepackage{overpic}
\usepackage{color}
\usepackage[square,comma,sort&compress,numbers]{natbib}
\usepackage{latexsym}
\usepackage{mathrsfs}
\usepackage{amsfonts}
\usepackage{hyperref}
\usepackage{amsthm}
\usepackage{appendix}
\usepackage{pdfsync}

\theoremstyle{definition}
\newtheorem{theorem}{Theorem}[section]
\newtheorem{lemma}[theorem]{Lemma}

\newtheorem{corollary}[theorem]{Corollary}
\newtheorem{definition}[theorem]{Definition}

\newtheorem{remark}[theorem]{Remark}
\newtheorem*{theorem*}{Theorem}
\def\qed{\hfill{Q.E.D.}\smallskip}

\begin{document}

\title{\bf A discrete uniformization theorem for decorated piecewise hyperbolic metrics on surfaces}
\author{Xu Xu, Chao Zheng}

\date{\today}

\address{School of Mathematics and Statistics, Wuhan University, Wuhan, 430072, P.R.China}
 \email{xuxu2@whu.edu.cn}

\address{School of Mathematics and Statistics, Wuhan University, Wuhan 430072, P.R. China}
\email{czheng@whu.edu.cn}



\begin{abstract}
In this paper, we study a natural discretization of the smooth Gaussian curvature on surfaces.
A discrete uniformization theorem is established for this discrete Gaussian curvature.
We further investigate the prescribing combinatorial curvature problem for a parametrization of this discrete Gaussian curvature,
which is called the combinatorial $\alpha$-curvature.
To find decorated piecewise hyperbolic metrics with prescribed combinatorial $\alpha$-curvatures,
we introduce the combinatorial $\alpha$-Ricci flow for decorated piecewise hyperbolic metrics.
To handle the potential singularities along the combinatorial $\alpha$-Ricci flow,
we do surgery along the flow by edge flipping under the weighted Delaunay condition.
Then we prove the longtime existence and convergence of the combinatorial $\alpha$-Ricci flow with surgery.
As an application of the combinatorial $\alpha$-Ricci flow with surgery, we give the existence of decorated piecewise hyperbolic metrics with prescribed combinatorial $\alpha$-curvatures.
We further introduce the combinatorial $\alpha$-Calabi flow with surgery and study its longtime behavior.
%These combinatorial curvature flows provide effective algorithms for finding decorated piecewise hyperbolic metrics with prescribed combinatorial $\alpha$-curvatures.
\end{abstract}

\maketitle

\section{Introduction}

This paper is a continuation of \cite{XZ}, in which we use variational principles with constraints to establish a discrete uniformization theorem for the combinatorial $\alpha$-curvature of decorated piecewise Euclidean metrics on surfaces.
In this paper, we use the combinatorial $\alpha$-curvature flows to establish the counterpart results for the combinatorial $\alpha$-curvature of decorated piecewise hyperbolic metrics on surfaces.

\subsection{Combinatorial curvature and discrete uniformization theorem}
The smooth Gaussian curvature at a point $p$ in a Riemann surface can be defined to be
\begin{equation*}
R(p)=\lim_{r\rightarrow 0}\frac{12}{\pi r^4}\big(\pi r^2-A(r)\big),
\end{equation*}
where $A(r)$ is the area of the geodesic disk of radius $r$ at $p$.
Applying this definition to a vertex $i$ of a piecewise hyperbolic surface gives a natural approximation of the smooth Gaussian curvature (up to a constant)
$\mathcal{R}_i=K_i/\sinh^2\frac{r_i}{2}$,
where $K_i=2\pi-\theta_i$ is the angle defect at $i$ and $\theta_i$ is the cone angle at $i$.
For the convenience of calculations, we use the following form
\begin{equation}\label{Eq: curvature R}
R_i=\frac{K_i}{\tanh^2\frac{r_i}{2}}.
\end{equation}
Note that $\tanh\frac{r}{2}\sim \sinh\frac{r}{2}\sim \frac{r}{2}$ as $r\rightarrow 0$,
$R_i$ is also an approximation of the smooth Gaussian curvature on surfaces.
We call $R_i$ the discrete Gaussian curvature or combinatorial curvature.
One of the main aims of this paper is to study the discrete uniformization theorem for the combinatorial curvature $R_i$.

Suppose $S$ is a connected closed surface and $V$ is a finite non-empty subset of $S$ with $|V|=N$.
The pair $(S,V)$ is called a marked surface and the points in $V$ are called vertices of the marked surface.
A piecewise hyperbolic (PH for short) metric $d$ on the marked surface $(S,V)$ is a hyperbolic cone metric with the conic singularities contained in $V$.
The marked surface $(S,V)$ endowed with a PH metric $d$ is called a PH surface, denoted by $(S,V,d)$.
A decoration $r$ on a PH surface $(S,V,d)$ is a choice of hyperbolic circle of radius $r_i$ at each vertex $i\in V$.
These circles in the decoration are called vertex-circles.
The pair $(d,r)$ is called a decorated PH metric on the marked surface $(S,V)$.
In this paper, we focus on the case that each pair of vertex-circles is separated.

\begin{theorem}\label{Thm: existence 1}
Let $(d,r)$ be a decorated PH metric on a marked surface $(S,V)$ with Euler characteristic $\chi(S)<0$ and let $\overline{R}: V\rightarrow(-\infty,0]$ be a given function defined on the vertices.
There exists a unique decorated PH metric $(\widetilde{d},\widetilde{r})$ discrete conformal equivalent to $(d,r)$ on $(S,V)$ with the prescribed combinatorial curvature $\overline{R}$.
\end{theorem}

For the definition of discrete conformal equivalence, please refer to Definition \ref{Def: DCE} for a marked surface with a fixed triangulation and to Definition \ref{Def: GDCE} for a marked surface with variable triangulations.
As a corollary of Theorem \ref{Thm: existence 1},
we have the following discrete uniformization theorem for the combinatorial curvature $R$ on decorated PH surfaces.
\begin{corollary}\label{Cor: existence 1}
For any decorated PH metric $(d,r)$ on a marked surface $(S,V)$ with Euler characteristic $\chi(S)<0$,
there exists a unique decorated PH metric $(\widetilde{d},\widetilde{r})$ discrete conformal equivalent to $(d,r)$ with non-positive constant combinatorial curvature $R$.
\end{corollary}

The combinatorial curvature $R$ in (\ref{Eq: curvature R}) was first introduced by Ge-Xu \cite{GX2} for Thurston's circle packing metrics on surfaces.
Since then, there are lots of works on different discretizations of the smooth Gaussian curvature on surfaces.
Following \cite{GX2}, we further introduce the following parameterized combinatorial curvature for decorated PH metrics on surfaces
\begin{equation}\label{Eq: curvature R_alpha}
R_{\alpha, i}=\frac{K_i}{\tanh^{-\alpha}\frac{r_i}{2}},
\end{equation}
where $\alpha\in \mathbb{R}$ is a constant.
If $\alpha=-2$, then $R_{-2, i}$ is the combinatorial curvature $R_i$ in (\ref{Eq: curvature R}).
We call $R_{\alpha}$ the combinatorial $\alpha$-curvature.

\begin{theorem}\label{Thm: main 1}
Suppose $(S,V)$ is a marked surface with a decorated PH metric $(d,r)$, $\alpha\in \mathbb{R}$ is a constant
and $\overline{R}: V\rightarrow(-\infty,2\pi)$ is a given function defined on the vertices.
If one of the following conditions is satisfied
\begin{description}
\item[(i)] $\alpha<0,\ \chi(S)<0,\ \overline{R}\leq 0$;
\item[(ii)] $\alpha>0,\ \overline{R}>0$;
\item[(iii)] (\cite{BL}) $\alpha=0$, $\sum^N_{i=1} \overline{R}_{i}>2\pi \chi(S)$.
\end{description}
Then there exists a unique decorated PH metric $(\widetilde{d},\widetilde{r})$ discrete conformal equivalent to $(d,r)$ with the prescribed combinatorial $\alpha$-curvature $\overline{R}$. 
\end{theorem}

Theorem \ref{Thm: main 1} is a generalization of Theorem \ref{Thm: existence 1}.
Specially, if $\alpha=-2$, then the case (i)
in Theorem \ref{Thm: main 1} is reduced to Theorem \ref{Thm: existence 1}.
By the relationship of the combinatorial $\alpha$-curvature $R_\alpha$ and the angle defect $K$ in (\ref{Eq: curvature R_alpha}),
the case (iii) in Theorem \ref{Thm: main 1} is covered by Bobenko-Lutz's work \cite{BL}.
In the following, we just prove the cases (i) and (ii) of Theorem \ref{Thm: main 1}.

As a corollary of Theorem \ref{Thm: main 1},
we have the following result on the existence of decorated PH metrics with constant combinatorial $\alpha$-curvatures.

\begin{corollary}
Suppose $(S,V)$ is a marked surface with a decorated PH metric $(d,r)$ and $\alpha\in \mathbb{R}$ is a constant.
\begin{description}
\item[(i)]
If $\alpha<0$ and $\chi(S)<0$, there exists a unique decorated PH metric $(\widetilde{d},\widetilde{r})$ discrete conformal equivalent to $(d,r)$ with non-positive constant combinatorial $\alpha$-curvature.
\item[(ii)]
If $\alpha>0$, there exists a unique decorated PH metric $(\widetilde{d},\widetilde{r})$ discrete conformal equivalent to $(d,r)$ with positive constant combinatorial $\alpha$-curvature.
\item[(iii)]
(\cite{BL}) If $\alpha=0$, for any constant $\overline{R}\in (\frac{2\pi\chi(S)}{N}, 2\pi)$,
there exists a unique decorated PH metric $(\widetilde{d},\widetilde{r})$ discrete conformal equivalent to $(d,r)$ with the combinatorial $\alpha$-curvature $\overline{R}$.
\end{description}
\end{corollary}

\begin{remark}
The discrete uniformization theorems for discrete conformal structures on polyhedral surfaces have been studied extensively.
Gu-Luo-Sun-Wu \cite{Gu1} and Gu-Guo-Luo-Sun-Wu \cite{Gu2} gave the discrete uniformization theorems for the vertex scaling of piecewise Euclidean (PE for short) metrics and PH metrics on surfaces respectively.
Using the combinatorial $\alpha$-curvature flows,
Xu \cite{Xu3} gave a parameterized discrete uniformization theorem for the vertex scaling of PE metrics,
which generalizes Gu-Luo-Sun-Wu's result \cite{Gu1}.
Following \cite{Xu3}, we \cite{XZ TAMS} also gave a parameterized discrete uniformization theorem for the vertex scaling of PH metrics,
which generalizes Gu-Guo-Luo-Sun-Wu's result \cite{Gu2}.
Motivated by Kou\v{r}imsk\'{a}'s work \cite{Kourimska},
we \cite{XZ CVPDE} prove the discrete Kazdan-Warner type theorem by variational principles with constraints, 
which generalizes Xu's result \cite{Xu3}.
Recently, Bobenko-Lutz \cite{BL,BL2} introduced decorated PE and PH metrics on surfaces and proved the corresponding discrete uniformization theorems.
Motivated by Bobenko-Lutz's work \cite{BL2}, we \cite{XZ0,XZ} proved some discrete uniformization theorems for decorated PE metrics for some modified combinatorial curvatures.
In fact, we prove the discrete Kazdan-Warner type theorem in \cite{XZ}, 
which generalizes Bobenko-Lutz's result \cite{BL2}.
\end{remark}

\subsection{Combinatorial $\alpha$-curvature flows}
Since Chow-Luo's pineering work \cite{Chow-Luo} on combinatorial Ricci flow for Thurston's circle packings,
combinatorial curvature flows have been an effective approach for finding polyhedral metrics with prescribed combinatorial curvatures.
There are lots of research activities on the combinatorial curvature flows for different discrete conformal structures on surfaces.
Following \cite{Chow-Luo,Ge1,Ge2}, we \cite{XZ2} introduced the combinatorial Ricci flow with surgery and the combinatorial Calabi flow with surgery for hyperbolic inversive distance circle packings on surfaces.
In this paper, we introduce the combinatorial $\alpha$-Ricci flow with surgery and the combinatorial $\alpha$-Calabi flow with surgery for decorated PH metrics on surfaces,
which is a parameterized generalization of the combinatorial curvature flows with surgery in \cite{XZ2}.

Set
\begin{equation}\label{Eq: h,r}
h_i=-\ln\tanh\frac{r_i}{2}
\end{equation}
for all $i\in V$.
The function $h: V\rightarrow \mathbb{R}_{>0}$ is called a discrete conformal factor in the following.
The combinatorial $\alpha$-curvature (\ref{Eq: curvature R_alpha}) can be written as
\begin{equation}\label{Eq: curvature R_alpha 2}
R_{\alpha, i}=\frac{K_i}{e^{\alpha h_i}}.
\end{equation}

Let $\mathcal{T}={(V,E,F)}$ be a triangulation of $(S, V)$, where $V,E,F$ are the sets of vertices, edges and faces respectively.
A triangulation $\mathcal{T}$ on a PH surface $(S,V, d)$ is a geodesic triangulation if the edges are geodesics in the PH metric $d$.
Denote a vertex, an edge and a face in the triangulation $\mathcal{T}$ by $i,\{ij\},\{ijk\}$ respectively.
A PH metric $d$ on a marked surface $(S,V)$ with a fixed geodesic triangulation $\mathcal{T}$ defines a map $l: E\rightarrow \mathbb{R}_{>0}$ such that $l_{ij}, l_{ik}, l_{jk}$ satisfy the triangle inequalities for any face $\{ijk\}\in F$.
Conversely, given a map $l: E\rightarrow \mathbb{R}_{>0}$ such that $l_{ij}, l_{ik}, l_{jk}$ satisfy the triangle inequalities for any face $\{ijk\}\in F$, every face in $F$ can be embedded in the hyperbolic space as a hyperbolic triangle with the edge length given by $l: E\rightarrow \mathbb{R}_{>0}$.
By gluing these hyperbolic triangles isometrically along the edges in pair,
one can obtain PH metrics on a triangulated surface $(S,V,\mathcal{T})$,
which give rise to PH metrics on $(S, V)$.
Therefore, we usually use $l: E\rightarrow \mathbb{R}_{>0}$ to denote a PH metric and use $(l,r)$ to denote a decorated PH metric on $(S,V)$ if the triangulation $\mathcal{T}$ is fixed.
%we usually use $l: E\rightarrow \mathbb{R}_{>0}$ to denote a PH metric on $(S,V)$ if the triangulation $\mathcal{T}$ is fixed.
%In the following, we use $l: E\rightarrow \mathbb{R}_{>0}$ to denote a PH metric and use $(l,r)$ to denote a decorated PH metric on a triangulated surface $(S,V,\mathcal{T})$.

\begin{definition}
Suppose $\mathcal{T}$ is a triangulation of $(S,V)$, $\alpha\in \mathbb{R}$ is a constant and $\overline{R}: V\rightarrow\mathbb{R}$ is a given function defined on the vertices.
The combinatorial $\alpha$-Ricci flow for decorated PH metrics on $(S,V,\mathcal{T})$ is defined to be
\begin{eqnarray}\label{Eq: CRF}
\begin{cases}
\frac{dh_i}{dt}=R_{\alpha,i}-\overline{R}_i,\\
h_i(0)=h_0.
\end{cases}
\end{eqnarray}
The combinatorial Calabi flow for decorated PH metrics on $(S,V,\mathcal{T})$ is defined to be
\begin{eqnarray}\label{Eq: CCF}
\begin{cases}
\frac{dh_i}{dt}
=\Delta^\mathcal{T}_\alpha(\overline{R}-R_\alpha)_i,\\
h_i(0)=h_0,
\end{cases}
\end{eqnarray}
where the discrete $\alpha$-Laplace operator $\Delta^\mathcal{T}_\alpha$ is defined by
\begin{equation}\label{Eq: alpha laplace}
\Delta^\mathcal{T}_\alpha f_i=\frac{1}{e^{\alpha h_i}}\sum_{j\in V}\frac{\partial K_i}{\partial h_j}f_j
\end{equation}
for any function $f: V\rightarrow \mathbb{R}$.
\end{definition}

Along the combinatorial $\alpha$-Ricci flow (\ref{Eq: CRF}) and the combinatorial $\alpha$-Calabi flow (\ref{Eq: CCF}),
singularities may develop,
including that some hyperbolic triangles degenerate and the discrete conformal factors $h$ tend to infinity.
These are supported by the discussions in \cite{GJ1, Xu2, Xu 21a, Xu 21b}.
To handle the potential singularities along these combinatorial $\alpha$-curvature flows,
we do surgery on the flows by edge flipping under the weighted Delaunay condition,
the idea of which comes from Bobenko-Lutz's recent work \cite{BL}.

Recall the following definition of weighted Delaunay triangulation.
For any hyperbolic triangle $\{ijk\}\in F$ with three vertex-circles attached to the vertices $i,j,k$,
there is a unique hyperbolic circle $C_{ijk}$ simultaneously orthogonal to the three vertex-circles.
Denote $\alpha_{ij}^k$ as the interior intersection angle of the hyperbolic circle $C_{ijk}$ and the edge $\{ij\}$.
For any adjacent hyperbolic triangles $\{ijk\}$ and $\{ijl\}$ sharing a common edge $\{ij\}$,
the edge $\{ij\}$ is \textit{weighted Delaunay} in the decorated PH metric $(l,r)$ if
$
\alpha_{ij}^k+\alpha_{ij}^l\leq \pi.
$
And the triangulation $\mathcal{T}$ is weighted Delaunay in $(l,r)$ if every edge in the triangulation $\mathcal{T}$ is weighted Delaunay.
There are other equivalent definitions for the weighted Delaunay triangulation using the signed distance of the center of $C_{ijk}$ to the edges.
Please refer to  \cite{CLXZ, Glickenstein DCG,Glickenstein JDG,Glickenstein preprint,GT} and others.

Along these combinatorial $\alpha$-curvature flows on $(S,V,\mathcal{T})$,
if $\mathcal{T}$ is weighted Delaunay in $(l(t),r(t))$ for $t\in [0,T]$ and not weighted Delaunay in $(l(t),r(t))$ for $t\in (T,T+\epsilon),\ \epsilon>0$,
there exists an edge $\{ij\}\in E$ such that $\alpha_{ij}^k+\alpha_{ij}^l\leq \pi$ for $t\in [0,T]$ and $\alpha_{ij}^k+\alpha_{ij}^l> \pi$ for $t\in (T,T+\epsilon)$.
Then we replace the triangulation $\mathcal{T}$ by a new triangulation $\mathcal{T}'$ at the time $t=T$ via replacing two hyperbolic triangles $\{ijk\}$ and $\{ijl\}$ adjacent to $\{ij\}$ by two new hyperbolic triangles $\{ikl\}$ and $\{jkl\}$.
This procedure is called a \textbf{surgery\  by\  flipping} on the triangulation $\mathcal{T}$,
which is also an isometry of $(S,V)$ with a decorated PH metric $(l(T),r(T))$.
After the surgery by flipping at the time $t=T$, we run these combinatorial $\alpha$-curvature flows on $(S,V,\mathcal{T}')$ with the initial decorated PH metric $(l(T),r(T))$.
Whenever the weighted Delaunay condition is not satisfied along these combinatorial $\alpha$-curvature flows,
we do surgery on these combinatorial $\alpha$-curvature flows by flipping.

We have the following result on the longtime existence and convergence for the solutions of the combinatorial $\alpha$-Ricci flow with surgery and the combinatorial $\alpha$-Calabi flow with surgery for decorated PH metrics on $(S,V)$.

\begin{theorem}\label{Thm: main 2}
Suppose $(S,V)$ is a marked surface with a decorated PH metric $(d,r)$, $\alpha\in \mathbb{R}$ is a constant and $\overline{R}: V\rightarrow (-\infty,2\pi)$ is a given function defined on the vertices.
If one of the conditions (i)(ii)(iii) in Theorem \ref{Thm: main 1} is satisfied,
then the solutions of the combinatorial $\alpha$-Ricci flow with surgery and the combinatorial $\alpha$-Calabi flow with surgery exist for all time and converge exponentially fast to  $\overline{h}$ such that the decorated PH metric $(d(\overline{h}),r(\overline{h}))$ discrete conformal equivalent to $(d,r)$ has the combinatorial $\alpha$-curvature $\overline{R}$.
\end{theorem}

\subsection{Organization of the paper}
The paper is organized as follows.
In Section \ref{Sec: FT}, we study the combinatorial $\alpha$-curvature flows on a marked surface with a fixed triangulation.
In Section \ref{Sec: VT}, we allow the triangulation of the marked surface to be changed by edge flipping and prove the global rigidity of the combinatorial $\alpha$-curvatures.
In Section \ref{Sec: CRFS}, we prove Theorem \ref{Thm: main 1} and the case of the combinatorial $\alpha$-Ricci flow with surgery in Theorem \ref{Thm: main 2}.
In Section \ref{Sec: CCFS}, we prove the case of the combinatorial $\alpha$-Calabi flow with surgery in Theorem \ref{Thm: main 2}.

\section{Combinatorial curvature flows for fixed triangulations}
\label{Sec: FT}

In this section, we first recall the definition of discrete conformal equivalence of decorated PH metrics on a triangulated surface.
Then we give the decomposition of the Jacobian $(\frac{\partial K}{\partial h})$,
which is useful in Section \ref{Sec: VT}.
In the end, we give the local convergence of the combinatorial $\alpha$-Ricci flow (\ref{Eq: CRF}) and the combinatorial $\alpha$-Calabi flow (\ref{Eq: CCF}).

\subsection{Discrete conformal equivalence}

For a marked surface $(S,V)$ with a fixed triangulation $\mathcal{T}$,
Bobenko-Lutz \cite{BL} introduced the following definition of discrete conformal equivalence for decorated PH metrics.

\begin{definition}[\cite{BL}, Proposition 2.7]
\label{Def: DCE}
Two decorated PH metrics $(l,r)$ and $(\widetilde{l},\widetilde{r})$ on a triangulated surface $(S,V,\mathcal{T})$ are discrete conformal equivalent
if and only if there exists a discrete conformal factor $u\in \mathbb{R}^N$ such that
\begin{equation}\label{Eq: DCE1}
\sinh \widetilde{r}_i=e^{u_i}\sinh r_i
\end{equation}
for $i\in V$ and
\begin{equation}\label{Eq: DCE2}
\begin{aligned}
\cosh \widetilde{l}_{ij}
=&e^{u_i+u_j}(\cosh l_{ij}-\cosh r_i\cosh r_j)\\
&+\sqrt{(1+e^{2u_i}\sinh^2 r_i)(1+e^{2u_j}\sinh^2 r_j)}
\end{aligned}
\end{equation}
for all $\{ij\}\in E$.
\end{definition}

\begin{remark}
The inversive distance of two hyperbolic circles $C_i$ and $C_j$ in the hyperbolic plane is defined as
\begin{equation}\label{Eq: inversive distance}
I_{ij}=\frac{\cosh l_{ij}-\cosh r_i\cosh r_j}{\sinh r_i\sinh r_j},
\end{equation}
where $l_{ij}$ is the hyperbolic distance of the centers of two hyperbolic circles $C_i, C_j$ and $r_i$, $r_j$ are the radii of $C_i, C_j$ respectively.
Note that the inversive distance is invariant under the M\"{o}bius transformations \cite{BL,Coxeter}.
Denote the inversive distance of two vertex-circles in $(l,r)$ and $(\widetilde{l},\widetilde{r})$ as $I$ and $\widetilde{I}$ respectively.
It is direct to check that if $I=\widetilde{I}$,
then $(l,r)$ and $(\widetilde{l},\widetilde{r})$ are discrete conformal equivalent in the sense of Definition \ref{Def: DCE}.
Conversely, if they are discrete conformal equivalent,
it is shown \cite{BL} that $I=\widetilde{I}$.
Since each pair of vertex-circles is required to be separated, it is easy to see that $I>1$ by (\ref{Eq: inversive distance}).
Therefore, the discrete conformal equivalent decorated PH metrics on triangulated surfaces in Definition \ref{Def: DCE} are exactly the separated hyperbolic inversive distance circle packings introduced by Bowers-Stephenson \cite{BS}.
\end{remark}

Suppose $(S,V,\mathcal{T})$ is a triangulated surface with inversive distance $I>1$.
In order that the edge lengths $l_{ij},l_{jk},l_{ki}$ defined by (\ref{Eq: inversive distance}) for $\{ijk\}\in F$ satisfy the triangle inequalities,
there are some restrictions on the radii $r_i, r_j, r_k$.
Denote the admissible space of the radius vectors $({r_i,r_j,r_k})\in \mathbb{R}^3_{>0}$ for a  hyperbolic triangle $\{ijk\}\in F$ by
\begin{equation*}
\mathcal{R}^{\mathcal{T}}_{ijk}=\{(r_i,r_j,r_k)\in \mathbb{R}^3_{>0}| l_{rs}+l_{rt}>l_{ts}, \{r,s,t\}=\{i,j,k\} \}.
\end{equation*}
By the map in (\ref{Eq: h,r}), we denote the image of $\mathcal{R}^{\mathcal{T}}_{ijk}$ under $h$ by $\Omega^{\mathcal{T}}_{ijk}$.
Set
\begin{equation}\label{Eq: AS2}
\Omega^{\mathcal{T}}
=\bigcap_{\{ijk\}\in F}\Omega^{\mathcal{T}}_{ijk}
\end{equation}
to be the admissible space of discrete conformal factors $h$ on $(S,V,\mathcal{T})$.

The following result obtained in \cite{Guo, L3, Xu1,Xu2} characterizes the admissible space $\Omega^{\mathcal{T}}_{ijk}$ and the Jacobian of the angle defect $K$ with respect to $h$.
Recall that the angle defect at $i\in V$ is $K_i=2\pi-\theta_i$,
where $\theta_i$ is the cone angle at $i$.
Denote $\theta^i_{jk}$ as the inner angle of the hyperbolic triangle $\{ijk\}$ at $i$.
Then the angle defect at $i\in V$ can be written as
$K_i=2\pi-\sum_{\{ijk\}\in F}\theta^i_{jk}$,
where the summation is taken over all hyperbolic triangles with $i$ as a vertex.

\begin{lemma}\label{Lem: matrix negative}
Suppose $(S,V,\mathcal{T})$ is a triangulated surface with a decorated PH metric $(l,r)$.
\begin{description}
\item[(i)]
The admissible space $\Omega^{\mathcal{T}}_{ijk}$ is a non-empty simply connected open set whose boundary is analytic.
This implies that the admissible space $\Omega^{\mathcal{T}}$ is a connected open subset of $\mathbb{R}^N_{>0}$.
\item[(ii)]
The matrix $\frac{\partial (\theta^{i}_{jk}, \theta^{j}_{ik}, \theta^{k}_{ij})}{\partial (h_i, h_j, h_k)}$ is symmetric and positive definite on $\Omega^{\mathcal{T}}_{ijk}$.
As a result, the matrix $L=\frac{\partial (K_1,..., K_N)}{\partial(h_1,...,h_N)}$ is symmetric and negative definite on $\Omega^\mathcal{T}$.
\end{description}
\end{lemma}

In \cite{BL}, Bobenko-Lutz gave the following decomposition for the matrix $L=(\frac{\partial K}{\partial h})$.

\begin{lemma}\label{Lem: L decomposition}
The negative definite matrix $L$ could be decomposed to be
$$L=L_A+L_B,$$
where $L_A=\text{diag}\{A_1,...,A_N\}$ is a diagonal matrix with
\begin{equation}\label{Eq: A i}
A_i=\frac{\partial}{\partial h_i}\text{Area}(S)
=\sum_{j\sim i}(-w_{ij})(\cosh l_{ij}-1)
\end{equation}
and $L_B$ is a symmetric matrix defined by
\begin{eqnarray}\label{Eq: L B}
(L_B)_{ij}=
\begin{cases}
-\sum_{k\sim i}w_{ik},  &{j=i},\\
w_{ij}, &{j\sim i},\\
0, &{\text{otherwise}}
\end{cases}
\end{eqnarray}
with
\begin{equation*}
w_{ij}
=-(\frac{\partial \theta_{jk}^i}{\partial h_j}
+\frac{\partial \theta_{jl}^i}{\partial h_j}).
\end{equation*}
Here $\{ijk\}$ and $\{ijl\}$ are two adjacent hyperbolic triangles sharing an edge $\{ij\}$.
Furthermore, if the triangulation $\mathcal{T}$ is weighted Delaunay in $(l,r)$,
then $w_{ij}\geq0$, $\sum_{j\sim i}w_{ij}>0$ and $A_i<0$.
This implies $L_B$ is negative semi-definite and $L_A$ is negative definite.
\end{lemma}

One can refer to \cite{XZ2} for a proof of Lemma \ref{Lem: L decomposition}.

\subsection{Local convergence of combinatorial $\alpha$-curvature flows}

The combinatorial $\alpha$-Ricci flow (\ref{Eq: CRF}) and the combinatorial $\alpha$-Calabi flow (\ref{Eq: CCF}) are ODE systems with smooth coefficients.
Hence, the solutions of these combinatorial curvature flows always exist locally around the initial time $t=0$.
We further have the following result on the longtime existence and convergence for the solutions of these combinatorial $\alpha$-curvature flows.

\begin{theorem}\label{Thm: local convergence}
Suppose $(S,V,\mathcal{T})$ is a triangulated surface with a decorated PH metric $(l,r)$, $\alpha\in \mathbb{R}$ is a constant and $\overline{R}: V\rightarrow\mathbb{R}$ is a given function defined on the vertices.
\begin{description}
\item[(i)]
If the solution of the combinatorial $\alpha$-Ricci flow (\ref{Eq: CRF}) or the combinatorial $\alpha$-Calabi flow (\ref{Eq: CCF}) converges,
then there exists a discrete conformal factor $\overline{h}\in \Omega^{\mathcal{T}}$ such that the combinatorial $\alpha$-curvature of the decorated PH metric $(l(\overline{h}),r(\overline{h}))$ is  $\overline{R}$.
\item[(ii)]
Suppose that there exists a discrete conformal factor $\overline{h}\in \Omega^{\mathcal{T}}$ such that the decorated PH metric $(l(\overline{h}),r(\overline{h}))$ has the combinatorial $\alpha$-curvature $\overline{R}$ and $\alpha \overline{R}\geq0$,
then there is a constant $\delta>0$ such that if $||R_\alpha(h(0))-R_\alpha(\overline{h})||<\delta$, the solution of the combinatorial $\alpha$-Ricci flow (\ref{Eq: CRF}) and the combinatorial $\alpha$-Calabi flow (\ref{Eq: CCF}) exists for all time and converges exponentially fast to $\overline{h}$ respectively.
\end{description}
\end{theorem}
\proof
\noindent\textbf{(i):}
Suppose the solution $h(t)$ of the combinatorial $\alpha$-Ricci flow (\ref{Eq: CRF}) converges to $\overline{h}$ as $t\rightarrow +\infty$,
then by the $C^1$-smoothness of $R_\alpha$,
we have $R_\alpha(\overline{h})=\lim_{t\rightarrow +\infty}R_\alpha(h(t))$.
Furthermore, there exists a sequence $t_n\in(n,n+1)$ such that
\begin{equation*}
h_i(n+1)-h_i(n)=h'_i(t_n)
=R_{\alpha,i}(h(t_n))-\overline{R}_i
\rightarrow 0,\ \text{as}\  n\rightarrow +\infty.
\end{equation*}
This implies $R_{\alpha,i}(\overline{h})
=\lim_{n\rightarrow +\infty}R_{\alpha,i}(h(t_n))
=\overline{R}_i$ for all $i\in V$
and $\overline{h}$ is a discrete conformal factor with the combinatorial $\alpha$-curvature $\overline{R}$.

Similar arguments can be applied to the combinatorial $\alpha$-Calabi flow (\ref{Eq: CCF}).
There exists a sequence $t_n\in(n,n+1)$ such that
\begin{equation*}
h_i(n+1)-h_i(n)=h'_i(t_n)
=\Delta^\mathcal{T}_\alpha(\overline{R}-
R_{\alpha}(h(t_n)))_i\rightarrow 0,\  \text{as}\  n\rightarrow +\infty.
\end{equation*}
The equality (\ref{Eq: alpha laplace}) yields $\Delta^\mathcal{T}_\alpha=\Sigma^{-\alpha}L$,
where $\Sigma=\text{diag}\{e^{h_1},,...,e^{h_N}\}$.
Lemma \ref{Lem: matrix negative} shows that $L$ is negative definite and
thus $\Delta^\mathcal{T}_\alpha$ is non-degenerate.
Hence, $R_{\alpha,i}(\overline{h})
=\lim_{n\rightarrow +\infty}R_{\alpha,i}(h(t_n))
=\overline{R}_i$ for all $i\in V$
and $\overline{h}$ is a discrete conformal factor with the combinatorial $\alpha$-curvature $\overline{R}$.

\noindent\textbf{(ii):}
Suppose there exists a discrete conformal factor $\overline{h}\in \Omega^{\mathcal{T}}$ such that $R_{\alpha}(\overline{h})=\overline{R}$.
For the combinatorial $\alpha$-Ricci flow (\ref{Eq: CRF}), set $\Gamma(h)=R_{\alpha}-\overline{R}$.
Direct calculations give
\begin{equation*}
D\Gamma|_{h=\overline{h}}
=\Sigma^{-\alpha}L-\alpha \Lambda =\Sigma^{-\frac{\alpha}{2}}
(\Sigma^{-\frac{\alpha}{2}}L\Sigma^{-\frac{\alpha}{2}}
-\alpha \Lambda)\Sigma^{\frac{\alpha}{2}},
\end{equation*}
where  $\Lambda=\text{diag}\{\overline{R}_1,...,\overline{R}_N\}$.
Combining Lemma \ref{Lem: matrix negative} and $\alpha \overline{R}\geq0$,  $D\Gamma|_{h=\overline{h}}$ has $N$ negative eigenvalues.
This implies that $\overline{h}$ is a local attractor of the combinatorial $\alpha$-Ricci flow (\ref{Eq: CRF}).
Then the conclusion follows from Lyapunov Stability Theorem (\cite{Pontryagin}, Chapter 5).

Similarly, for the combinatorial $\alpha$-Calabi flow (\ref{Eq: CCF}), set $\Gamma(h)=\Delta^\mathcal{T}_\alpha(\overline{R}
-R_\alpha)$.
Direct calculations give
\begin{equation*}
\begin{aligned}
D\Gamma|_{h=\overline{h}}
&=-\Sigma^{-\alpha}L\Sigma^{-\alpha}L
+\alpha\Sigma^{-\alpha}L\Lambda \\
&=-\Sigma^{-\alpha}(-L)\Sigma^{-\alpha}(-L)
-\alpha\Sigma^{-\alpha}(-L)\Lambda \\
&=-\Sigma^{-\frac{\alpha}{2}}
\left[\Sigma^{-\frac{\alpha}{2}}(-L)\Sigma^{-\alpha}
(-L)\Sigma^{-\frac{\alpha}{2}}+\alpha
\Sigma^{-\frac{\alpha}{2}} (-L)\Sigma^{-\frac{\alpha}{2}} \Lambda \right]\Sigma^{\frac{\alpha}{2}}\\
&=-\Sigma^{-\frac{\alpha}{2}}\left(Q^2+\alpha Q \Lambda\right)\Sigma^{\frac{\alpha}{2}}\\
&=-\Sigma^{-\frac{\alpha}{2}}Q^{\frac{1}{2}}
\left[Q^2+Q^{\frac{1}{2}}(\alpha \Lambda)Q^{\frac{1}{2}}\right]Q^{-\frac{1}{2}}
\Sigma^{\frac{\alpha}{2}},
\end{aligned}
\end{equation*}
where $Q=\Sigma^{-\frac{\alpha}{2}} (-L)\Sigma^{-\frac{\alpha}{2}}$.
By Lemma \ref{Lem: matrix negative}, $Q$ is a symmetric and positive definite matrix.
Combining this with $\alpha \overline{R}\geq0$,  $D\Gamma|_{h=\overline{h}}$ has $N$ negative eigenvalues.
This implies that $\overline{h}$ is a local attractor of the combinatorial $\alpha$-Calabi flow (\ref{Eq: CCF}).
The conclusion follows from Lyapunov Stability Theorem (\cite{Pontryagin}, Chapter 5).
\qed

Theorem \ref{Thm: local convergence} gives the longtime existence and convergence for the solutions of the combinatorial $\alpha$-curvature flows on a triangulated surface for initial decorated PH metrics with small initial energy.
However, for general initial decorated PH metrics, the combinatorial $\alpha$-curvature flows may develop singularities.
To handle the possible singularities along the combinatorial $\alpha$-curvature flows,
we need to do surgery by flipping on the combinatorial $\alpha$-curvature flows.

\section{Rigidity of the combinatorial $\alpha$-curvatures}
\label{Sec: VT}

To analyze the longtime behavior of the combinatorial $\alpha$-Ricci flow with surgery and the combinatorial $\alpha$-Calabi flow with surgery,
we need the discrete conformal theory for decorated PH metrics recently established by Bobenko-Lutz \cite{BL}.
In this section, we first briefly recall Bobenko-Lutz's discrete conformal theory for decorated PH metrics.
Then we extend the definition of combinatorial $\alpha$-curvatures to decorated PH surfaces and prove the global rigidity of the combinatorial $\alpha$-curvatures on decorated PH surfaces.

\subsection{Bobenko-Lutz's discrete conformal theory}
The following definition of discrete conformal equivalence of decorated PH metrics generalizes Definition \ref{Def: DCE},
which allows the triangulation of the marked surface to be changed under the weighted Delaunay condition.

\begin{definition}[\cite{BL}, Proposition 4.9]\label{Def: GDCE}
Let $(d,r)$ and $(\widetilde{d},\widetilde{r})$ be two decorated PH metrics on the marked surface $(S,V)$.
They are discrete conformal equivalent if and only if there is a finite sequence of triangulated decorated PH surfaces
$(\mathcal{T}^0,l^0,r^0),...,(\mathcal{T}^m,l^m,r^m)$ such that
\begin{description}
\item[(i)]
the decorated PH metrics of $(\mathcal{T}^0,l^0,r^0)$ and $(\mathcal{T}^m,l^m,r^m)$ are $(d,r)$ and $(\widetilde{d},\widetilde{r})$ respectively,
\item[(ii)]
each $\mathcal{T}^n$ is a weighted Delaunay triangulation of the decorated PH surface $(\mathcal{T}^n,l^n,r^n)$,
\item[(iii)]
if $\mathcal{T}^n=\mathcal{T}^{n+1}$, then there is a discrete conformal factor $u\in \mathbb{R}^N$ such that $(\mathcal{T}^n,l^n,r^n)$ and $(\mathcal{T}^{n+1},l^{n+1},r^{n+1})$ are related by (\ref{Eq: DCE1}) and (\ref{Eq: DCE2}),
\item[(iv)]
if $\mathcal{T}^n\neq\mathcal{T}^{n+1}$, then $\mathcal{T}^n$ and $\mathcal{T}^{n+1}$ are two different weighted Delaunay triangulations of the same decorated PH surface.
\end{description}
\end{definition}

%Definition \ref{Def: GDCE} gives an equivalence for the decorated PH metrics on a marked surface.
%The equivalence class of a decorated PH metric $(d,r)$ on $(S,V)$ is denoted by $\mathcal{D}(d,r)$.
Using this definition of discrete conformal equivalence, Bobenko-Lutz \cite{BL} proved the following important theorem for decorated PH metrics.

\begin{theorem}[\cite{BL}, Theorem A]\label{Thm: BL}
Let $(d,r)$ be a decorated PH metric on the marked surface $(S,V)$.
There exists a unique decorated PH metric $(\widetilde{d},\widetilde{r})$ discrete conformal equivalent to $(d,r)$ realizing $\overline{K}: V\rightarrow (-\infty, 2\pi)$ if and only if $\overline{K}$ satisfies the discrete Gauss-Bonnet condition $\sum_{i=1}^N K_i>2\pi \chi(S).$
\end{theorem}

Submitting (\ref{Eq: DCE1}) into (\ref{Eq: h,r}) gives the relationship of $h$ and $u$, i.e.,
\begin{equation}\label{Eq: h,u}
\coth h_i(u)=\cosh r_i(u)
=\sqrt{1+\sinh^2 r_i(u)}
=\sqrt{1+e^{2u_i}\sinh^2 r_i(0)},
\end{equation}
where $r_i(0)$ is the initial value.
Hence, the map from $h$ to $u$ is a bijection.
Denote $\mathcal{D}(d,r)$ as the equivalence class of a decorated PH metric $(d,r)$ on $(S,V)$.
For simplicity, for $(\widetilde{d},\widetilde{r})\in \mathcal{D}(d,r)$,
we sometimes denote it by $(d(h),r(h))$ for some $h\in \mathbb{R}^N_{>0}$.
Set
\begin{equation*}
\mathcal{C}_\mathcal{T}(d,r)
=\{h\in \mathbb{R}^N_{>0} |\ \mathcal{T}\ \text{is a weighted Delaunay triangulation of}\ (S,V,d(h),r(h))\}.
\end{equation*}
Note that $\mathcal{C}_\mathcal{T}(d,r)$ is a subset of the admissible space $\Omega^{\mathcal{T}}$ in (\ref{Eq: AS2}).
Furthermore, Bobenko-Lutz \cite{BL} gave the following result.

\begin{lemma}[\cite{BL}]
\label{Lem: finite decomposition}
The set
$$J=\{\mathcal{T}| \mathcal{C}_{\mathcal{T}}(d,r)\ \text{has non-empty interior in}\ \mathbb{R}_{>0}^N\}$$
is a finite set, $\mathbb{R}_{>0}^N=\cup_{\mathcal{T}_i\in J}\mathcal{C}_{\mathcal{T}_i}(d,r)$ and
each $\mathcal{C}_{\mathcal{T}_i}(d,r)$ is homeomorphic to $\mathbb{R}^N_{>0}$.
\end{lemma}

Combining Lemma \ref{Lem: matrix negative} (ii) and Lemma \ref{Lem: finite decomposition},
the function $\mathcal{W}_\mathcal{T}:\mathcal{C}_{\mathcal{T}}(d,r)\rightarrow \mathbb{R}$ defined by
\begin{equation*}
\mathcal{W}_\mathcal{T}(h)
=\int^h\sum_{i=1}^NK_idh_i
\end{equation*}
is a well-defined locally strictly concave function  on $\mathcal{C}_\mathcal{T}(d,r)$.
To extend the function $\mathcal{W}_\mathcal{T}$, we need the following lemma obtained from \cite{L3}.
\begin{definition}[\cite{L3}, Definition 2.3]
A differential 1-form $w=\sum_{i=1}^n a_i(x)dx^i$ in an open subset $U\subset \mathbb{R}^N$ is said to be continuous if each $a_i(x)$ is continuous on $U$.
A continuous differential 1-form $w$ is called closed if $\int_{\partial \tau}w=0$ for each triangle $\tau\subset U$.
\end{definition}

\begin{lemma}[\cite{L3}, Proposition 2.5]
\label{Lem: L3}
Suppose $X$ is an open convex subset of $\mathbb{R}^N$ and $A\subseteq X$ is an open subset bounded by a $(n-1)$-dimensional $C^1$-smooth submanifold in $X$.
If $w=\sum_{i=1}^n a_i(x)dx^i$ is a continuous closed 1-form on $X$ such that $F(x)=\int_a^xw$ is locally convex in $A$ and $X-\overline{A}$,
where $\overline{A}$ is the closure of $A$ in $X$,
then $F(x)$ is convex in $X$.
\end{lemma}

By Lemma \ref{Lem: finite decomposition}, one can extend the combinatorial map $K$ defined on $\mathcal{C}_{\mathcal{T}}(d,r)$ to be
\begin{equation}\label{Eq: K extended}
\begin{aligned}
\mathbf{K}: \mathbb{R}_{>0}^N&\rightarrow (-\infty,2\pi)^N\\
h&\mapsto K(h)
\end{aligned}
\end{equation}
defined on $\mathbb{R}_{>0}^N$,
which is independent of the choice of the weighted Delaunay triangulations on $(S,V)$.
Hence, using Lemma \ref{Lem: L3}, one can extend the function $\mathcal{W}_\mathcal{T}$ to be the following function
\begin{equation}\label{Eq: W}
\mathcal{W}(h)
=\int^h\sum_{i=1}^N\mathbf{K}_idh_i
\end{equation}
defined on $\mathbb{R}_{>0}^N$.
Note that $\mathcal{W}(h)$ is a $C^1$-smooth concave function on $\mathbb{R}_{>0}^N$.
Bobenko-Lutz \cite{BL} further proved the following result.

\begin{lemma}(\cite{BL}, Proposition 4.11)
\label{Lem: W}
The function $\mathcal{W}(h)$ is a $C^2$-smooth strictly concave function defined on $\mathbb{R}_{>0}^N$.
\end{lemma}

\begin{remark}
It should be mentioned that the function $\mathcal{W}(h)$ is equivalent to the discrete Hilbert-Einstein functional $\mathcal{H}_{\Sigma_g,\Theta}(h)$ introduced by Bobenko-Lutz \cite{BL} up to a constant.
\end{remark}

Lemma \ref{Lem: W} implies that $\mathbf{K}$ defined on $\mathbb{R}_{>0}^N$ is a $C^1$-extension of the combinatorial curvature $K$ defined on $\mathcal{C}_\mathcal{T}(d,r)$.
Note that the discrete $\alpha$-Laplace operator $\Delta^\mathcal{T}_\alpha$ is independent of the choice of the weighted Delaunay triangulations of a decorated PH metric, i.e., it is intrinsic.
Hence, the discrete Laplace operator $\Delta_\alpha^\mathcal{T}$ could be extend to the following operator $\Delta_\alpha$ defined on $\mathbb{R}_{>0}^N$,
which is continuous and piecewise smooth on $\mathbb{R}_{>0}^N$ as a matrix-valued function of $h$.

\begin{definition}\label{Def: Delta alpha}
Let $(d,r)$ be a decorated PH metric on the marked surface $(S,V)$.
The discrete $\alpha$-Laplace operator $\Delta_\alpha$ is defined to be the map
\begin{equation*}
\begin{aligned}
\Delta_\alpha :\  &\mathbb{R}_{>0}^N\rightarrow \mathbb{R}^N\\
&f\mapsto \Delta_\alpha^\mathcal{T} f,
\end{aligned}
\end{equation*}
with
\begin{equation}\label{Eq: F67}
\Delta_\alpha f_i
=\frac{1}{e^{\alpha h_i}}\sum_{j\in V}\frac{\partial \mathbf{K}_i}{\partial h_j}f_j
=\frac{1}{e^{\alpha h_i}}(\mathbf{L}f)_i
\end{equation}
for any $f : V \rightarrow \mathbb{R}$, where $\mathbf{L}_{ij}=\frac{\partial \mathbf{K}_i}{\partial h_j}$ is an extension of $L_{ij}=\frac{\partial K_i}{\partial h_j}$.
\end{definition}

\begin{remark}
By Lemma \ref{Lem: L decomposition}, we have
\begin{equation}\label{Eq: Delta alpha f_i}
\begin{aligned}
\Delta_\alpha f_i
&=\frac{1}{e^{\alpha h_i}}(\mathbf{L}f)_i\\
&=\frac{1}{e^{\alpha h_i}}\mathbf{L}_{ii}f_i
+\frac{1}{e^{\alpha h_i}}\sum_{j\neq i,\ j\in V}\mathbf{L}_{ij}f_j\\
&=\frac{1}{e^{\alpha h_i}}\big(\mathbf{A}_i+(\mathbf{L}_B)_{ii}\big)f_i
+\frac{1}{e^{\alpha h_i}}\sum_{j\neq i,\ j\in V}(\mathbf{L}_B)_{ij}f_j\\
&=\sum_{j\sim i}\frac{\mathbf{w}_{ij}}{e^{\alpha h_i}}(f_j-f_i)+\frac{\mathbf{A}_i}{e^{\alpha h_i}}f_i
\end{aligned}
\end{equation}
for any $f : V \rightarrow \mathbb{R}$,
where $\mathbf{w}_{ij}=\mathbf{L}_{ij}$ for $j\sim i$ is an extension of $w_{ij}$ and $\mathbf{A}_i$ is an extension of $A_i$.
\end{remark}

\subsection{Rigidity of the combinatorial $\alpha$-curvatures}

The following combinatorial $\alpha$-curvature $\mathbf{R}_{\alpha}$ for decorated PH metrics on a marked surface
is an extension of the combinatorial $\alpha$-curvature $R_\alpha$ in (\ref{Eq: curvature R_alpha 2}) on a triangulated surface.

\begin{definition}\label{Def: R alpha extended}
Suppose $(S,V)$ is a marked surface with a decorated PH metric $(d,r)$ and $\alpha\in \mathbb{R}$ is a constant.
The combinatorial $\alpha$-curvature at the vertex $i\in V$ for $(\widetilde{d},\widetilde{r})\in \mathcal{D}(d,r)$ is defined to be
\begin{equation}\label{Eq: R alpha extended}
\mathbf{R}_{\alpha, i}
=\frac{\mathbf{K}_i}{e^{\alpha h_i}},
\end{equation}
where $\mathbf{K}_i$ is the combinatorial curvature at $i$ defined by (\ref{Eq: K extended}).
\end{definition}

A basic problem on the combinatorial $\alpha$-curvature is to understand the relationships between the decorated PH metrics and the combinatorial $\alpha$-curvatures.
The following theorem shows the global rigidity of decorated PH metrics with respect to the combinatorial $\alpha$-curvatures,
which corresponds to the rigidity part of Theorem \ref{Thm: main 1}.

\begin{theorem}\label{Thm: rigidity}
Suppose $(S,V)$ is a marked surface with a decorated PH metric $(d,r)$, $\alpha\in \mathbb{R}$ is a constant and $\overline{\mathbf{R}}: V\rightarrow\mathbb{R}$ is a given function defined on the vertices.
If $\alpha\overline{\mathbf{R}}\geq 0$, then there exists at most one discrete conformal factor $h^*\in \mathbb{R}^N_{>0}$ such that the decorated PH metric $(d(h^*),r(h^*))\in \mathcal{D}(d,r)$ has the combinatorial $\alpha$-curvature $\overline{\mathbf{R}}$.
\end{theorem}
\proof
By Lemma \ref{Lem: W}, the following function
\begin{equation*}
\mathcal{W}_\alpha(h)
=\mathcal{W}(h)-\int^h
\sum_{i=1}^N\overline{\mathbf{R}}_ie^{\alpha h_i}dh_i
=\int^h\sum_{i=1}^N\mathbf{K}_idh_i-\int^h
\sum_{i=1}^N\overline{\mathbf{R}}_ie^{\alpha h_i}dh_i
\end{equation*}
is a $C^2$-smooth function on $\mathbb{R}^N_{>0}$.
By direct calculations, we have
\begin{equation*}
\nabla_{h_i}\mathcal{W}_\alpha(h)
=\mathbf{K}_i-\overline{\mathbf{R}}_ie^{\alpha h_i}
=(\mathbf{R}_{\alpha, i}-\overline{\mathbf{R}}_i)e^{\alpha h_i},
\end{equation*}
where (\ref{Eq: R alpha extended}) is used in the second equality.
Hence for $h^*\in \mathbb{R}^N_{>0}$, the decorated PH metric $(d(h^*),r(h^*))\in \mathcal{D}(d,r)$ has the combinatorial $\alpha$-curvature $\overline{\mathbf{R}}$ if and only if $\nabla_{h_i}\mathcal{W}(h^*)=0$ for all $i\in V$.
Moreover,
\begin{equation*}%\label{Eq: F2}
\mathrm{Hess}_h\ \mathcal{W}_\alpha
=\mathbf{L}-\alpha \left(
 \begin{array}{ccc}
  \overline{\mathbf{R}}_1e^{\alpha h_1}  &  &\\
     & \ddots &  \\
      & & \overline{\mathbf{R}}_Ne^{\alpha h_N} \\
  \end{array}
 \right).
\end{equation*}
Lemma \ref{Lem: matrix negative} shows that $\mathbf{L}$ is negative definite, i.e., $\mathbf{L}<0$.
Combining with the assumption that $\alpha\overline{\mathbf{R}}\geq 0$,
we have $\mathrm{Hess}_h\ \mathcal{W}_\alpha<0$.
This implies $\mathcal{W}_\alpha$ is strictly concave on $\mathbb{R}^N_{>0}$.
The conclusion follows from the following classical result from calculus.

\noindent\textbf{Lemma:}
If $f:\Omega \rightarrow \mathbb{R}$ is a $C^1$-smooth  strictly concave function on an open convex set $\Omega \subset \mathbb{R}^N$,
then its gradient $\nabla f:\Omega \rightarrow \mathbb{R}^N$ is injective.
\qed

\section{Combinatorial $\alpha$-Ricci flow with surgery}
\label{Sec: CRFS}

In this section, we first introduce the combinatorial $\alpha$-Ricci flow with surgery on a decorated PH surface.
Then we prove that the convergence of the combinatorial $\alpha$-Ricci flow with surgery is equivalent to the existence of decorated PH metrics with the prescribed combinatorial $\alpha$-curvature.
In the end, we use the maximal principle to prove the existence of decorated PH metrics with the prescribed combinatorial $\alpha$-curvatures, i.e., Theorem \ref{Thm: main 1}.

\begin{definition}\label{Def: CRFS}
Suppose $(S,V)$ is a marked surface with a decorated PH metric $(d,r)$, $\alpha\in \mathbb{R}$ is a constant and $\overline{\mathbf{R}}: V\rightarrow\mathbb{R}$ is a given function defined on the vertices.
The combinatorial $\alpha$-Ricci flow with surgery is defined to be
\begin{eqnarray}\label{Eq: CRFS}
\begin{cases}
\frac{dh_i}{dt}=\mathbf{R}_{\alpha, i}-\overline{\mathbf{R}}_i,\\
h_i(0)=h_0.
\end{cases}
\end{eqnarray}
\end{definition}

\subsection{An equivalent theorem}

\begin{theorem}\label{Thm: equivalent 1}
Suppose $(S,V)$ is a marked surface with a decorated PH metric $(d,r)$, $\alpha\in \mathbb{R}$ is a constant and $\overline{\mathbf{R}}: V\rightarrow \mathbb{R}$ is a given function defined on the vertices.
\begin{description}
\item[(i)]
If the solution of the combinatorial $\alpha$-Ricci flow with surgery (\ref{Eq: CRFS}) converges,
then there exists a decorated PH metric in $\mathcal{D}(d,r)$ with the combinatorial $\alpha$-curvature $\overline{\mathbf{R}}$.
\item[(ii)]
If $\overline{\mathbf{R}}: V\rightarrow(-\infty,2\pi)$ satisfies $\alpha\overline{\mathbf{R}}\geq0$ and
there exists a decorated PH metric in $\mathcal{D}(d,r)$ with the combinatorial $\alpha$-curvature $\overline{\mathbf{R}}$,
then the solution of the combinatorial $\alpha$-Ricci flow with surgery (\ref{Eq: CRFS}) exists for all time and converges exponentially fast for any initial $h(0)\in \mathbb{R}_{>0}^N$.
\end{description}
\end{theorem}
\proof
\noindent\textbf{(i):}
The proof is almost the same as that of Theorem \ref{Thm: local convergence} (i), so we omit it here.

\noindent\textbf{(ii):}
Combining the assumptions and Theorem \ref{Thm: rigidity},
there exists a unique discrete conformal factor $\overline{h}$ such that $\mathbf{R}_\alpha(\overline{h})=\overline{\mathbf{R}}$.
Set
\begin{equation}\label{Eq: W alpha}
\mathcal{W}_\alpha(h)
=\mathcal{W}(h)-\int^h_{\overline{h}}
\sum_{i=1}^N\overline{\mathbf{R}}_ie^{\alpha h_i}dh_i
=\int^h_{\overline{h}}\sum_{i=1}^N\mathbf{K}_idh_i -\int^h_{\overline{h}}
\sum_{i=1}^N\overline{\mathbf{R}}_ie^{\alpha h_i}dh_i
\end{equation}
By the proof of Theorem \ref{Thm: rigidity},
$\mathcal{W}_\alpha(h)$ is a $C^2$-smooth strictly concave function on $\mathbb{R}^N_{>0}$.
Moreover, $\mathcal{W}_\alpha(\overline{h})=0,\ \nabla \mathcal{W}_\alpha(\overline{h})=0$.
This implies $0=\mathcal{W}(\overline{h})\geq \mathcal{W}(h)$ and $\lim_{||h||\rightarrow +\infty}\mathcal{W}_\alpha(h)=-\infty$ by the concavity of $\mathcal{W}_\alpha(h)$.
Hence, $\mathcal{W}(u)$ is proper.

Along the combinatorial $\alpha$-Ricci flow with surgery (\ref{Eq: CRFS}), we have
\begin{equation*}
\frac{d\mathcal{W}_\alpha(h(t))}{dt}
=\sum^N_{i=1}\frac{\partial \mathcal{W}_\alpha}{\partial h_i}\frac{dh_i}{dt}
=\sum^N_{i=1}(\mathbf{R}_{\alpha,i}
-\overline{\mathbf{R}}_i)e^{\alpha h_i}
(\mathbf{R}_{\alpha,i}-\overline{\mathbf{R}}_i)
=\sum^N_{i=1}(\mathbf{R}_{\alpha,i}-\overline{\mathbf{R}}_i)^2
e^{\alpha h_i}\geq 0.
\end{equation*}
This implies $\mathcal{W}_\alpha(h(0))\leq \mathcal{W}_\alpha(h(t))\leq 0$.
Combining this with the properness of $\mathcal{W}_\alpha(h)$,
the solution $h(t)$ of the combinatorial $\alpha$-Ricci flow with surgery (\ref{Eq: CRFS}) is uniformly bounded in $\mathbb{R}^N$.
Note that $h\in \mathbb{R}^N_{>0}$,
then $h(t)$ is uniformly bounded from above in $\mathbb{R}^N_{>0}$.
The following Lemma \ref{Lem: below 1} shows that $h(t)$ is uniformly bounded from below in $\mathbb{R}^N_{>0}$.
Then the solution $h(t)$ of the combinatorial $\alpha$-Ricci flow with surgery (\ref{Eq: CRFS}) lies in a compact subset of $\mathbb{R}^N_{>0}$.
This implies the solution of the combinatorial $\alpha$-Ricci flow with surgery (\ref{Eq: CRFS}) exists for all time and $\lim_{t\rightarrow +\infty}\mathcal{W}_\alpha(h(t))$ exists.

By the mean value theorem,
there exists a sequence $t_n\in(n,n+1)$ such that as $n\rightarrow +\infty$,
\begin{equation*}
\begin{aligned}
\mathcal{W}_\alpha(h(n+1))-\mathcal{W}_\alpha(h(n))
&=(\mathcal{W}_\alpha(h(t)))'|_{t_n}
=\nabla \mathcal{W}_\alpha\cdot\frac{dh_i}{dt}|_{t_n}\\
&=\sum^N_{i=1}(\mathbf{R}_{\alpha,i}
-\overline{\mathbf{R}}_i)^2e^{\alpha h_i}|_{t_n}
\rightarrow 0.
\end{aligned}
\end{equation*}
Then $\lim_{n\rightarrow +\infty}\mathbf{R}_{\alpha,i}(h(t_n))
=\overline{\mathbf{R}}_i
=\mathbf{R}_{\alpha,i}(\overline{h})$ for all $i\in V$.
By $\{h(t)\}\subset\subset \mathbb{R}^N_{>0}$,
there exists $h^*\in \mathbb{R}^N_{>0}$ and a convergent subsequence of $\{h(t_n)\}$, still denoted as $\{h(t_n)\}$ for simplicity, such that $\lim_{n\rightarrow \infty}h(t_n)=h^*$.
This implies
$\mathbf{R}_\alpha(h^*)=\lim_{n\rightarrow +\infty}\mathbf{R}_\alpha(h(t_n))
=\mathbf{R}_\alpha(\overline{h})$.
Then $h^*=\overline{h}$ by Theorem \ref{Thm: rigidity}.
Therefore, $\lim_{n\rightarrow \infty}h(t_n)=\overline{h}$.

Similar to the proof of Theorem \ref{Thm: local convergence} (ii),
set $\Gamma(h)=\mathbf{R}_{\alpha}-\overline{\mathbf{R}}$,
then $D\Gamma|_{h=\overline{h}}$ has $N$ negative eigenvalues.
This implies that $\overline{h}$ is a local attractor of the combinatorial $\alpha$-Ricci flow with surgery (\ref{Eq: CRFS}) .
Then the conclusion follows from Lyapunov Stability Theorem (\cite{Pontryagin}, Chapter 5).
\qed

\begin{lemma}\label{Lem: below 1}
Suppose $(S,V)$ is a marked surface with a decorated PH metric $(d,r)$, $\alpha\in \mathbb{R}$ is a constant and $\overline{\mathbf{R}}: V\rightarrow(-\infty,2\pi)$ is a given function defined on the vertices and satisfies $\alpha\overline{\mathbf{R}}\geq0$.
Along the combinatorial $\alpha$-Ricci flow with surgery (\ref{Eq: CRFS}),
$h_i(t)$ is uniformly bounded from below in $\mathbb{R}_{>0}$ for all $i\in V$.
\end{lemma}

To prove Lemma \ref{Lem: below 1},
we need the following lemma,
which has been proved for different cases in almost the same manner.
One can refer to \cite{Xu 21a,BL} for the proof.

\begin{lemma}\label{Lem: BL3}
Let $(S,V,\mathcal{T})$ be a weighted Delaunay triangulated surface with a decorated PH metric $(l,r)$.
For any $\epsilon>0$, there exists a constant $C=C(\epsilon)>0$ such that if $h_i<C$, then $\theta^i_{jk}<\epsilon$ for any hyperbolic triangle $\{ijk\}\in F$ with $i$ as a vertex.
\end{lemma}

\noindent\textbf{Proof of Lemma \ref{Lem: below 1}:}
Lemma \ref{Lem: finite decomposition} implies the finiteness of the weighted Delaunay triangulations
along the combinatorial Ricci flow with surgery (\ref{Eq: CRFS}).
Hence, we just need to prove this lemma for some fixed weighted Delaunay triangulation $\mathcal{T}$ in $(d,r)$.

We prove this by contradiction.
Otherwise, there exists at least one vertex $i\in V$ such that $\lim_{t\rightarrow T}h_i(t)=0$ for $T\in(0,+\infty]$.
\begin{description}
\item[(i)]
If $\alpha<0$, then $\overline{R}\leq0$ by the assumption.
By Lemma \ref{Lem: BL3}, choosing $\epsilon_1=\frac{2\pi}{N}>0$,
there exists a positive constant $C_1=C_1(\epsilon_1)$ such that whenever $h_i<C_1$,
the angle $\theta^i_{jk}<\epsilon_1$.
Then
\begin{equation*}
K_i=2\pi-\sum_{\{ijk\}\in F}\theta^i_{jk}
>2\pi-\sum_{\{ijk\}\in F}\frac{2\pi}{N}>0.
\end{equation*}
As a result, $R_{\alpha,i}
=\frac{K_i}{e^{\alpha h_i}}>0\geq \overline{R}_i$.

\item[(ii)]
If $\alpha=0$, then $\overline{R}\in (-\infty,2\pi)$ by the assumption.
Choosing $\epsilon_2=\frac{1}{N}(2\pi-\overline{R}_i)>0$,
there exists a positive constant $C_2=C_2(\epsilon_2)$ such that whenever $h_i<C_2$,
the angle $\theta^i_{jk}<\epsilon_2$.
Then
$K_i>2\pi-\sum_{\{ijk\}\in F}\frac{1}{N}(2\pi-\overline{R}_i)>\overline{R}_i$.
Hence $R_{\alpha,i}
=\frac{K_i}{e^{\alpha h_i}}=K_i>\overline{R}_i$.

\item[(iii)]
If $\alpha>0$, then $\overline{R}\in [0,2\pi)$ by the assumption.
If $\overline{R}\equiv0$, then by the same arguments as that in the case (i), we have $R_{\alpha,i}>\overline{R}_i$ if $h_i<C_1$.
Suppose $\overline{R}\geq0$ and $\overline{R}\not\equiv0$.
Set $\overline{R}_{\max}=\max_{i\in V}\overline{R}_i$.
Then $0<\overline{R}_{\max}<2\pi$.
By Lemma \ref{Lem: BL3}, choosing a sufficient small $\epsilon_3$ with $0<\epsilon_3<\frac{1}{N}(2\pi-\overline{R}_{\max})$,
there exists a positive constant $C_3=C_3(\epsilon_3)$ such that whenever $h_i<C_3$,
the angle $\theta^i_{jk}<\epsilon_3$.
Then
\begin{equation*}
K_i=2\pi-\sum_{\{ijk\}\in F}\theta^i_{jk}
>2\pi-N\epsilon_3>\overline{R}_{\max}.
\end{equation*}
Set $C_4=\min\{\frac{1}{\alpha}\ln \frac{2\pi-N\epsilon_3}{\overline{R}_{\max}},C_3\}$.
If $h_i<C_4$, we have
\begin{equation*}
R_{\alpha,i}
=\frac{K_i}{e^{\alpha h_i}}
>\frac{K_i}{e^{\alpha C_4}}
\geq\frac{K_i}{\frac{2\pi-N\epsilon_3}{\overline{R}_{\max}}}
>\overline{R}_{\max}
\geq\overline{R}_i.
\end{equation*}
\end{description}
In summary, there exists a constant $C_5=\min\{C_1,C_2,C_4\}>0$ such that if $h_i<C_5$, then $R_{\alpha,i}> \overline{R}_i$.

Choose a time $t_0\in (0,T)$ such that $h_i(t_0)<C_5$, this can be done because $\lim_{t\rightarrow T}h_i(t)=0$.
Set $a=\inf\{t<t_0|h_i(s)<C_5, \forall s\in (t,t_0]\}$, then $h_i(a)=C_5$.
By the arguments above, we have $\frac{d h_i}{dt}=R_{\alpha,i}-\overline{R}_i>0$ on $(a,t_0]$.
This implies $h_i(t_0)>h_i(a)=C_5$, which contradicts $h_i(t_0)<C_5$.
Hence, $h_i(t)$ is uniformly bounded from below in $\mathbb{R}_{>0}$ for any $i\in V$ along the combinatorial $\alpha$-Ricci flow with surgery (\ref{Eq: CRFS}).
\qed

\begin{remark}
Specially, if $\alpha= 0$, then $\mathbf{R}_{\alpha,i}=\mathbf{K}_i$.
In this case, if we further require that $\overline{\mathbf{R}}_i$ satisfies the discrete Gauss-Bonnet condition $\sum_{i=1}^N \overline{\mathbf{R}}_i>2\pi \chi(S)$,
then by Theorem \ref{Thm: BL}, there exists a decorated PH metric with the combinatorial curvature $\overline{\mathbf{R}}$.
Using Theorem \ref{Thm: equivalent 1} (ii), we see that
the solution of the combinatorial $\alpha$-Ricci flow with surgery (\ref{Eq: CRFS}) exists for all time and converges exponentially fast for any initial $h(0)\in \mathbb{R}_{>0}^N$.
This has been proved in \cite{XZ2}.
\end{remark}

\subsection{Existence of decorated PH metrics with the prescribed combinatorial $\alpha$-curvatures}
Set
\begin{equation}\label{Eq: M alpha}
M_{\alpha,i}=\mathbf{R}_{\alpha,i}
-\overline{\mathbf{R}}_i.
\end{equation}
Submitting (\ref{Eq: M alpha}) into (\ref{Eq: CRFS}) gives
\begin{equation}\label{Eq: CRFS2}
\frac{dh_i}{dt}=M_{\alpha,i}.
\end{equation}
The following lemma gives the evolution equation for $M_{\alpha,i}$ along the combinatorial $\alpha$-Ricci flow with surgery (\ref{Eq: CRFS2}).

\begin{lemma}
Along the combinatorial $\alpha$-Ricci flow with surgery (\ref{Eq: CRFS2}),
$M_{\alpha,i}$ evolves according to
\begin{equation}\label{Eq: M alpha i}
\frac{dM_{\alpha,i}}{dt}
=\sum_{j\sim i}\frac{\mathbf{w}_{ij}}{e^{\alpha h_i}}(M_{\alpha,j}-M_{\alpha,i}) +\frac{\mathbf{A}_i}{e^{\alpha h_i}}M_{\alpha,i}-\alpha M_{\alpha,i}(M_{\alpha,i}+\overline{\mathbf{R}}_i),
\end{equation}
where $\mathbf{w}_{ij}\geq0$ and $\mathbf{A}_i<0$.
\end{lemma}
\proof
By direct calculations, we have
\begin{equation*}
\begin{aligned}
\frac{dM_{\alpha,i}}{dt}
=&\sum^N_{j=1}\frac{\partial \mathbf{R}_{\alpha,i}}{\partial h_j}\frac{dh_j}{dt}\\
=&\sum^N_{j=1}(\frac{1}{e^{\alpha h_i}}\frac{\partial \mathbf{K}_i}{\partial h_j}
-\frac{\mathbf{K}_i}{e^{\alpha h_i}}\alpha \delta_{ij})\cdot M_{\alpha,j}\\
=&\sum^N_{j=1}\frac{1}{e^{\alpha h_i}}\frac{\partial \mathbf{K}_i}{\partial h_j}M_{\alpha,j}
-\mathbf{R}_{\alpha,i}\alpha M_{\alpha,i}\\
=&(\Delta_\alpha M_\alpha)_i-\alpha M_{\alpha,i}(M_{\alpha,i}+\overline{\mathbf{R}}_i)\\
=&\sum_{j\sim i}\frac{\mathbf{w}_{ij}}{e^{\alpha h_i}}(M_{\alpha,j}-M_{\alpha,i}) +\frac{\mathbf{A}_i}{e^{\alpha h_i}}M_{\alpha,i}-\alpha M_{\alpha,i}(M_{\alpha,i}+\overline{\mathbf{R}}_i),
\end{aligned}
\end{equation*}
where the second line uses (\ref{Eq: R alpha extended}) and (\ref{Eq: CRFS2}),
the fourth line uses (\ref{Eq: F67}) and (\ref{Eq: M alpha}), and
the last line uses (\ref{Eq: Delta alpha f_i}).
Moreover, $\mathbf{w}_{ij}\geq0$, $\mathbf{A}_i<0$ follows from Lemma \ref{Lem: L decomposition}.
\qed

To analyze the longtime behavior of the combinatorial $\alpha$-Ricci flow with surgery (\ref{Eq: CRFS2}),
we need the following discrete maximum principle obtained in \cite{GX2}.

\begin{theorem}(\cite{GX2}, Theorem 3.5)
\label{Thm: MX}
Let $f:V\times [0,T)\rightarrow \mathbb{R}$ be a $C^1$ function such that
\begin{equation*}
\frac{\partial f_i}{\partial t}\geq \sum_{j\sim i}a_{ij}(f_j-f_i)+\Phi_i(f_i),\ \forall(v_i,t)\in V\times [0,T),
\end{equation*}
where $a_{ij}\geq 0$ and $\Phi_i: \mathbb{R}\rightarrow \mathbb{R}$ is a local Lipschitz function.
Suppose there exists $C_1\in \mathbb{R}$ such that $f_i(0)\geq C_1$ for all $i\in V$.
Let $\varphi$ be the solution to the associated ODE
\begin{eqnarray*}
\begin{cases}
\frac{d\varphi}{dt}=\Phi_i(\varphi),\\
\varphi(0)=C_1,
\end{cases}
\end{eqnarray*}
then $f_i(t)\geq \varphi(t)$ for all $(i,t)\in V\times [0,T)$ such that $\varphi(t)$ exists.

Similarly, suppose $f:V\times [0,T)\rightarrow \mathbb{R}$ be a $C^1$ function such that
\begin{equation*}
\frac{\partial f_i}{\partial t}\leq \sum_{j\sim i}a_{ij}(f_j-f_i)+\Phi_i(f_i),\ \forall(v_i,t)\in V\times [0,T).
\end{equation*}
Suppose there exists $C_2\in \mathbb{R}$ such that $f_i(0)\leq C_2$ for all $i\in V$.
Let $\psi$ be the solution to the associated ODE
\begin{eqnarray*}
\begin{cases}
\frac{d\psi}{dt}=\Phi_i(\psi),\\
\psi(0)=C_2,
\end{cases}
\end{eqnarray*}
then $f_i(t)\leq \psi(t)$ for all $(i,t)\in V\times [0,T)$ such that $\psi(t)$ exists.
\end{theorem}

\begin{remark}
As  mentioned in Remark 3 of \cite{GX2}, $\sum_{j\sim i}a_{ij}(f_j-f_i)$ in Theorem \ref{Thm: MX} is a generalization of the classical discrete Laplace operator $\Delta f_i$ in \cite{Chung}.
Specially, the symmetry condition $a_{ij}=a_{ji}$ is required for the classical discrete Laplace operator $\Delta f_i$ in \cite{Chung}, but not for $\sum_{j\sim i}a_{ij}(f_j-f_i)$ in Theorem \ref{Thm: MX}.
This is important for the applications in the following.
%Here we use the notation $\sum_{j\sim i}a_{ij}(f_j-f_i)$ rather than $\Delta f_i$.
%One reason is to avoid confusion with the $\alpha$-Laplace operator $\Delta_\alpha$ in Definition \ref{Def: Delta alpha},
%the others is that $\sum_{j\sim i}a_{ij}(f_j-f_i)$ is formally the same as $\sum_{j\sim i}\frac{\mathbf{w}_{ij}}{e^{\alpha h_i}}(M_{\alpha,j}-M_{\alpha,i})$.
\end{remark}

Applying Theorem \ref{Thm: MX} to the evolution equation (\ref{Eq: M alpha i}) gives the following result for $M_{\alpha}$.
The proof is direct and thus we omit it here.
Indeed, $0$ is the solution of the ODE.

\begin{corollary}\label{Cor: application 1}
Along the combinatorial $\alpha$-Ricci flow with surgery (\ref{Eq: CRFS2}),
if $M_{\alpha,i}(0)\geq 0$ for all $i\in V$, then $M_{\alpha,i}(t)\geq 0$ for all $i\in V$;
if $M_{\alpha,i}(0)\leq 0$ for all $i\in V$, then $M_{\alpha,i}(t)\leq 0$ for all $i\in V$.
\end{corollary}
%\proof
%For any three constants $a,b,c$ and $a\neq 0$, we consider the following ODE
%\begin{equation*}
%\frac{d\varphi}{dt}
%=a\varphi-b \varphi(\varphi+c).
%\end{equation*}
%The solution is
%\begin{equation*}
%\varphi(t)=
%\begin{cases}
%\frac{a-bc}{b+(a-bc)\mu e^{(bc-a)t}}, &{b\neq 0},\\
%\lambda e^{at}, &{b=0},
%\end{cases}
%\end{equation*}
%where $\mu,\lambda$ are two constants.
%If $\varphi(0)=0$,
%then $a-bc=0$ for $b\neq 0$ and $\lambda=0$ for $b=0$.
%This implies that $\varphi(t)\equiv 0$ for any $t$.
%\qed

For $\alpha<0$, we have the following result on the existence of decorated PH metrics with the prescribed combinatorial $\alpha$-curvature.

\begin{lemma}\label{Lem: case1}
Suppose $(S,V)$ is a marked surface with a decorated PH metric $(d,r)$, $\alpha<0$ is a constant and $\overline{\mathbf{R}}: V\rightarrow(-\infty,0]$ is a given function defined on the vertices.
If one of the following two conditions holds
\begin{description}
\item[(i)]
$\overline{\mathbf{R}}$ is a negative constant and there exists $h^\ast$ such that $\overline{\mathbf{R}}<\mathbf{R}_{\alpha,i}(h^\ast)<0$ for all $i\in V$,
\item[(ii)]
$\overline{\mathbf{R}}$ is a non-positive constant and there exists $h^\ast$ such that $\mathbf{R}_{\alpha}(h^\ast)$ is a negative constant,
\end{description}
then there exists a decorated PH metric in $\mathcal{D}(d,r)$ with the combinatorial $\alpha$-curvature $\overline{\mathbf{R}}$.
\end{lemma}
\proof
\noindent\textbf{(i):}
Take $h^*$ as the initial value of the combinatorial $\alpha$-Ricci flow with surgery (\ref{Eq: CRFS2}).
Set
\begin{equation*}
\mathbf{R}_{\alpha,\max}(h(0))=\max_{i\in V}\mathbf{R}_{\alpha,i}(h^\ast),\ M_{\alpha,\max}=\max_{i\in V}M_{\alpha,i}.
\end{equation*}
\begin{equation*}
\mathbf{R}_{\alpha,\min}(h(0))=\min_{i\in V}\mathbf{R}_{\alpha,i}(h^\ast),\  M_{\alpha,\min}=\min_{i\in V}M_{\alpha,i}.
\end{equation*}
The assumption implies that $\overline{\mathbf{R}}<\mathbf{R}_{\alpha,\min}(h(0))<0$.
Then $M_{\alpha,i}(0)>0$ for all $i\in V$.
By Corollary \ref{Cor: application 1}, we have
$M_{\alpha,i}(t)\geq 0$ for all $i\in V$.
This implies that 
\begin{equation*}
\frac{dh_i}{dt}=M_{\alpha,i}(t)\geq 0 
\end{equation*}
for all $i\in V$ by (\ref{Eq: CRFS2}).
Then $h_i(t)\geq h_i(0)$ for all $i\in V$.
Due to the fact that monotone bounded function must converge,
we just need to prove that $h_i(t)$ is uniformly bounded from above in $\mathbb{R}_{>0}$ for all $i\in V$.
Then the solution $h_i(t)$ of the combinatorial $\alpha$-Ricci flow with surgery (\ref{Eq: CRFS2}) with initial value $h(0)=h^*$ lies in a compact subset of $\mathbb{R}_{>0}^N$.
This implies that the solution of the combinatorial $\alpha$-Ricci flow with surgery (\ref{Eq: CRFS2}) with initial value $h(0)=h^*$ exists for all time and converges. 
By Theorem \ref{Thm: equivalent 1} (i),
there exists a decorated PH metric in $\mathcal{D}(d,r)$ with the combinatorial $\alpha$-curvature $\overline{\mathbf{R}}$.

Now we start to prove that $h_i(t)$ is uniformly bounded from above in $\mathbb{R}_{>0}$ for all $i\in V$.
Submitting $\mathbf{A}_i<0$ and $M_{\alpha,i}(t)\geq 0$ into (\ref{Eq: M alpha i}) gives
\begin{equation*}
\begin{aligned}
\frac{dM_{\alpha,i}}{dt}
&=\sum_{j\sim i}\frac{\mathbf{w}_{ij}}{e^{\alpha h_i}}(M_{\alpha,j}-M_{\alpha,i}) +\frac{\mathbf{A}_i}{e^{\alpha h_i}}M_{\alpha,i}-\alpha M_{\alpha,i}(M_{\alpha,i}+\overline{\mathbf{R}})\\
&\leq \sum_{j\sim i}\frac{\mathbf{w}_{ij}}{e^{\alpha h_i}}(M_{\alpha,j}-M_{\alpha,i})
-\alpha M_{\alpha,i}(M_{\alpha,i}+\overline{\mathbf{R}}).
\end{aligned}
\end{equation*}
Consider the following ODE
\begin{eqnarray*}
\begin{cases}
\frac{d \psi}{dt}=-\alpha \psi(\psi+\overline{\mathbf{R}}),\\
\psi(0)=M_{\alpha,\max}(0).
\end{cases}
\end{eqnarray*}
Note that $\overline{\mathbf{R}}<0$,
the solution of the ODE is
\begin{equation*}
\psi(t)
=\frac{\overline{\mathbf{R}}}
{-1+c_1e^{\alpha\overline{\mathbf{R}}t}},
\end{equation*}
where $c_1$ is a constant.
At $t=0$, $M_{\alpha,\max}(0)=\frac{\overline{\mathbf{R}}}
{-1+c_1}$.
Then $c_1=\frac{\overline{\mathbf{R}}}{M_{\alpha,\max}(0)}+1
=\frac{\mathbf{R}_{\alpha,\max}(h(0))}{M_{\alpha,\max}(0)}<0$.
By Theorem \ref{Thm: MX}, we have
\begin{equation*}
M_{\alpha,i}(t)
\leq \frac{\mathbf{\overline{R}}}
{-1+(\frac{\mathbf{\overline{R}}}{M_{\alpha,\max}(0)}+1)
e^{\alpha \mathbf{\overline{R}}t}}
\leq \frac{\overline{\mathbf{R}}M_{\alpha,\max}(0)}
{\mathbf{R}_{\alpha,\max}(h(0))}e^{-\alpha \mathbf{\overline{R}}t}
\end{equation*}
for all $i\in V$.
Set $c_2=\frac{\overline{\mathbf{R}}M_{\alpha,\max}(0)}
{\mathbf{R}_{\alpha,\max}(h(0))}>0$, then
$\frac{dh_i}{dt}=M_{\alpha,i}(t)\leq c_2e^{-\alpha \mathbf{\overline{R}}t}$ by (\ref{Eq: CRFS2}).
Hence,
\begin{equation*}
h_i(t)\leq h_i(0)+\frac{c_2}{\alpha \mathbf{\overline{R}}}(1-e^{-\alpha \mathbf{\overline{R}}t}).
\end{equation*}
Combining with $\alpha \mathbf{\overline{R}}>0$, $h_i(t)$ is uniformly bounded from above in $\mathbb{R}_{>0}$ for all $i\in V$.

\noindent\textbf{(ii):}
According to the assumptions, we see that the relationship of $\mathbf{R}_{\alpha}(h^*)$ and $\overline{\mathbf{R}}$ is
$$\overline{\mathbf{R}}=\mathbf{R}_{\alpha}(h^*)<0\quad \text{or} \quad
\overline{\mathbf{R}}<\mathbf{R}_{\alpha}(h^*)<0 \quad \text{or} \quad
\mathbf{R}_{\alpha}(h^*)<\overline{\mathbf{R}}\leq 0.$$
If $\overline{\mathbf{R}}=\mathbf{R}_{\alpha}(h^*)<0$, then $(d(h^*),r(h^*))$ is the decorated PH metric with the combinatorial $\alpha$-curvature $\overline{\mathbf{R}}$ and thus
the conclusion holds.
If $\overline{\mathbf{R}}<\mathbf{R}_{\alpha}(h^*)<0$, the conclusion follows from the case (i).
If $\mathbf{R}_{\alpha}(h^*)<\overline{\mathbf{R}}\leq0$, we take $h^*$ as the initial value of the combinatorial $\alpha$-Ricci flow with surgery (\ref{Eq: CRFS2}).
Then $M_{\alpha}(0)<0$ is a constant.
By Corollary \ref{Cor: application 1}, we have $M_{\alpha,i}(t)\leq 0$ for all $i\in V$.
This implies that 
\begin{equation}\label{Eq: F8}
\frac{dh_i}{dt}=M_{\alpha,i}(t)\leq 0 
\end{equation}
for all $i\in V$ by (\ref{Eq: CRFS2}).
Then $h_i(t)\leq h_i(0)$ for all $i\in V$.
Similar to the arguments in the case (i), 
we just need to prove that $h_i(t)$ is uniformly bounded from below in $\mathbb{R}_{>0}$ for all $i\in V$.

Submitting $\mathbf{A}_i<0$ and $M_{\alpha,i}(t)\leq0$ into (\ref{Eq: M alpha i}) gives
\begin{equation*}
\begin{aligned}
\frac{dM_{\alpha,i}}{dt}
&=\sum_{j\sim i}\frac{\mathbf{w}_{ij}}{e^{\alpha h_i}}(M_{\alpha,j}-M_{\alpha,i}) +\frac{\mathbf{A}_i}{e^{\alpha h_i}}M_{\alpha,i}-\alpha M_{\alpha,i}(M_{\alpha,i}+\overline{\mathbf{R}})\\
&\geq \sum_{j\sim i}\frac{\mathbf{w}_{ij}}{e^{\alpha h_i}}(M_{\alpha,j}-M_{\alpha,i})
-\alpha M_{\alpha,i}(M_{\alpha,i}+\overline{\mathbf{R}}).
\end{aligned}
\end{equation*}
Consider the following ODE
\begin{eqnarray*}
\begin{cases}
\frac{d \varphi}{dt}=-\alpha \varphi(\varphi+\overline{\mathbf{R}}),\\
\varphi(0)=M_{\alpha}(0).
\end{cases}
\end{eqnarray*}
If $\mathbf{\overline{R}}<0$, by similar calculations in the case (i), the solution of the ODE is
\begin{equation*}
\varphi(t)=\frac{\mathbf{\overline{R}}}
{-1+(\frac{\mathbf{\overline{R}}}{M_{\alpha}(0)}+1)
e^{\alpha \mathbf{\overline{R}}t}}.
\end{equation*}
If $\mathbf{\overline{R}}=0$, the solution of the ODE is
\begin{equation*}
\varphi(t)=\frac{M_{\alpha}(0)}{\alpha M_{\alpha}(0)t+1}.
\end{equation*}

Hence, by Theorem \ref{Thm: MX}, if $\mathbf{\overline{R}}<0$, then
\begin{equation*}
M_{\alpha,i}(t) \geq
\frac{\mathbf{\overline{R}}}
{-1+(\frac{\mathbf{\overline{R}}}{M_{\alpha}(0)}+1)
e^{\alpha \mathbf{\overline{R}}t}}
\geq M_{\alpha}(0)e^{-\alpha \mathbf{\overline{R}}t}.
\end{equation*}
By (\ref{Eq: CRFS2}), we have
$\frac{dh_i}{dt}=M_{\alpha,i}(t)\geq M_{\alpha}(0)e^{-\alpha\mathbf{\overline{R}}t}$.
This implies
\begin{equation}\label{Eq: F5}
h_i(t)\geq h_i(0)+\frac{M_{\alpha}(0)}{\alpha \mathbf{\overline{R}}}(1-e^{-\alpha \mathbf{\overline{R}}t}).
\end{equation}
Thus $h_i(t)\geq h_i(0)+\frac{M_{\alpha}(0)}{\alpha \mathbf{\overline{R}}}$ for all $(i,t)\in V\times [0,+\infty)$.
Since $M_{\alpha}(0)<0$ and $h_i>0$, this uniform lower bound may be non-valid.
Note that Lemma \ref{Lem: below 1} shows that $h_i(t)$ is uniformly bounded from below in $\mathbb{R}_{>0}$ for all $i\in V$, denoted by $C_6$.
If $h_i(0)+\frac{M_{\alpha}(0)}{\alpha \mathbf{\overline{R}}}\geq C_6$ for all $i\in V$, 
then $h_i(t)\in [h_i(0)+\frac{M_{\alpha}(0)}{\alpha \mathbf{\overline{R}}},h_i(0)]$ for all $i\in V$.
Combining with (\ref{Eq: F8}), the solution of the combinatorial $\alpha$-Ricci flow with surgery (\ref{Eq: CRFS2}) with initial value $h(0)=h^*$ exists for all time and converges.
The conclusion follows from Theorem \ref{Thm: equivalent 1} (i).
If $h_i(0)+\frac{M_{\alpha}(0)}{\alpha \mathbf{\overline{R}}}< C_6$ for all $i\in V$, 
then $h_i(t)\in [C_6,h_i(0)]$ for all $i\in V$.
Combining with (\ref{Eq: F8}), we see that
there exists $T_1$ such that $\frac{dh_i}{dt}\leq 0$ for all $(i,t)\in V\times [0,T_1)$ and $\frac{dh_i}{dt}\equiv0$ for all $(i,t)\in V\times [T_1,+\infty)$.
In this case, the solution of the combinatorial $\alpha$-Ricci flow with surgery (\ref{Eq: CRFS2}) with initial value $h(0)=h^*$ exists for all time and converges.
The conclusion follows from Theorem \ref{Thm: equivalent 1} (i).

By Theorem \ref{Thm: MX}, if $\mathbf{\overline{R}}=0$, then
\begin{equation*}
\frac{dh_i}{dt}=M_{\alpha,i}(t) \geq\frac{M_{\alpha}(0)}{\alpha M_{\alpha}(0)t+1}.
\end{equation*}
This implies
\begin{equation}\label{Eq: F6}
h_i(t)\geq h_i(0)+\frac{1}{\alpha}\ln (\alpha M_{\alpha}(0)t+1).
\end{equation}
Since $\alpha<0$, then $h_i(0)+\frac{1}{\alpha}\ln (\alpha M_{\alpha}(0)t+1)\rightarrow -\infty$ as $t\rightarrow +\infty$.
Then $h_i(t)\in [C_6,h_i(0)]$ for all $i\in V$.
Similar the above arguments,  there exists $T_2$ such that $\frac{dh_i}{dt}\leq 0$ for all $(i,t)\in V\times [0,T_2)$ and $\frac{dh_i}{dt}\equiv0$ for all $(i,t)\in V\times [T_2,+\infty)$.
In this case, the solution of the combinatorial $\alpha$-Ricci flow with surgery (\ref{Eq: CRFS2}) with initial value $h(0)=h^*$ exists for all time and converges.
The conclusion follows from Theorem \ref{Thm: equivalent 1} (i).
\qed

Similar to Lemma \ref{Lem: case1}, we have the following result for $\alpha>0$.

\begin{lemma}\label{Lem: case2}
Suppose $(S,V)$ is a marked surface with a decorated PH metric $(d,r)$, $\alpha>0$ is a constant and $\overline{\mathbf{R}}: V\rightarrow(0,2\pi)$ is a given positive constant function defined on the vertices.
If there exists $h^\ast$ such that one of the following two conditions holds
\begin{description}
\item[(i)]
$0<\overline{\mathbf{R}}<\mathbf{R}_{\alpha,i}(h^\ast)$ for all $i\in V$,
\item[(ii)]
$\mathbf{R}_{\alpha}(h^\ast)$ is a positive constant,
\end{description}
then there exists a decorated PH metric in $\mathcal{D}(d,r)$ with the combinatorial $\alpha$-curvature $\overline{\mathbf{R}}$.
\end{lemma}
\proof
\noindent\textbf{(i):}
Take $h^*$ as the initial value of the combinatorial $\alpha$-Ricci flow with surgery (\ref{Eq: CRFS2}).
The assumption implies that $0<\overline{\mathbf{R}}<\mathbf{R}_{\alpha,\min}(h(0))$.
Then $M_{\alpha,i}(0)>0$ for all $i\in V$.
By Corollary \ref{Cor: application 1}, we have
$M_{\alpha,i}(t)\geq 0$ for all $i\in V$.
Submitting $\mathbf{A}_i<0$ and $M_{\alpha,i}(t)\geq 0$ into (\ref{Eq: M alpha i}) gives
\begin{equation*}
\frac{dM_{\alpha,i}}{dt}
\leq \sum_{j\sim i}\frac{\mathbf{w}_{ij}}{e^{\alpha h_i}}(M_{\alpha,j}-M_{\alpha,i})
-\alpha M_{\alpha,i}(M_{\alpha,i}+\overline{\mathbf{R}}).
\end{equation*}
Note that $\overline{\mathbf{R}}>0$, then by Theorem \ref{Thm: MX}, we have
\begin{equation*}
0\leq M_{\alpha,i}(t)
\leq \frac{\mathbf{\overline{R}}}
{-1+(\frac{\mathbf{\overline{R}}}{M_{\alpha,\max}(0)}+1)
e^{\alpha \mathbf{\overline{R}}t}}
\leq M_{\alpha,\max}(0)e^{-\alpha \mathbf{\overline{R}}t}.
\end{equation*}
The rest proof is parallelling to that of Lemma \ref{Lem: case1} (i), so we omit it here.

\noindent\textbf{(ii):}
Because $\mathbf{R}_{\alpha}(h^*)$ and $\overline{\mathbf{R}}$ are both constants, then
$$0<\overline{\mathbf{R}}=\mathbf{R}_{\alpha}(h^*)\quad \text{or} \quad
0<\overline{\mathbf{R}}<\mathbf{R}_{\alpha}(h^*)
\quad \text{or} \quad
0<\mathbf{R}_{\alpha}(h^*)<\overline{\mathbf{R}}.$$
If $0<\overline{\mathbf{R}}=\mathbf{R}_{\alpha}(h^*)$, then $(d(h^*),r(h^*))$ is the decorated PH metric with the combinatorial $\alpha$-curvature $\overline{\mathbf{R}}$ and thus
the conclusion holds.
If $0<\overline{\mathbf{R}}<\mathbf{R}_{\alpha}(h^*)$, the conclusion follows from the case (i).
If $0<\mathbf{R}_{\alpha}(h^*)<\overline{\mathbf{R}}$, we take $h^*$ as the initial value of the combinatorial $\alpha$-Ricci flow with surgery (\ref{Eq: CRFS2}).
Then $M_{\alpha}(0)<0$ is a constant.
By Corollary \ref{Cor: application 1}, we have $M_{\alpha,i}(t)\leq 0$ for all $i\in V$.
Submitting $\mathbf{A}_i<0$ and $M_{\alpha,i}(t)\leq0$ into (\ref{Eq: M alpha i}) gives
\begin{equation*}
\frac{dM_{\alpha,i}}{dt}
\geq \sum_{j\sim i}\frac{\mathbf{w}_{ij}}{e^{\alpha h_i}}(M_{\alpha,j}-M_{\alpha,i})
-\alpha M_{\alpha,i}(M_{\alpha,i}+\overline{\mathbf{R}}).
\end{equation*}
Note that $\overline{\mathbf{R}}>0$, then by Theorem \ref{Thm: MX}, we have
\begin{equation*}
0\geq M_{\alpha,i}(t) \geq
\frac{\mathbf{\overline{R}}}
{-1+(\frac{\mathbf{\overline{R}}}{M_{\alpha}(0)}+1)
e^{\alpha \mathbf{\overline{R}}t}}
\geq \frac{\overline{\mathbf{R}}M_{\alpha}(0)}
{\mathbf{R}_{\alpha}(h(0))}e^{-\alpha \mathbf{\overline{R}}t}.
\end{equation*}
The rest proof is parallelling to that of Lemma \ref{Lem: case1} (ii), so we omit it here.
\qed

\begin{remark}
Note that we do not have a result under the condition that $\alpha>0,\ 0=\overline{\mathbf{R}}<\mathbf{R}_{\alpha,i}(h^*)$ for all $i\in V$.
The main reason is that we do not have a result paralleling to Lemma \ref{Lem: BL3} under the condition that $h_i$ is bounded from below for all $i\in V$.
Indeed, under the condition that $\alpha>0$, $0=\overline{\mathbf{R}}<\mathbf{R}_{\alpha,i}(h^*)$ for all $i\in V$,
we have $M_{\alpha,i}(0)>0$ for all $i\in V$.
By Corollary \ref{Cor: application 1}, we have
$M_{\alpha,i}(t)\geq 0$ for all $i\in V$.
Then
\begin{equation*}
\frac{dM_{\alpha,i}}{dt}
\leq \sum_{j\sim i}\frac{\mathbf{w}_{ij}}{e^{\alpha h_i}}(M_{\alpha,j}-M_{\alpha,i})
-\alpha M_{\alpha,i}(M_{\alpha,i}+\overline{\mathbf{R}}).
\end{equation*}
Note that $\overline{\mathbf{R}}=0$, by Theorem \ref{Thm: MX}, we have
\begin{equation*}
0\leq M_{\alpha,i}(t)=\frac{dh_i}{dt} \leq\frac{M_{\alpha}(0)}{\alpha M_{\alpha}(0)t+1}.
\end{equation*}
Then
\begin{equation*}
0\leq h_i(t)-h_i(0)\leq \frac{1}{\alpha}\ln (\alpha M_{\alpha}(0)t+1),
\end{equation*}
which can not ensure that $h_i(t)$ is uniformly bounded from above in $\mathbb{R}_{>0}$ for all $i\in V$ along the combinatorial $\alpha$-Ricci flow with surgery (\ref{Eq: CRFS2}).
\end{remark}

As a direct corollary of Lemma \ref{Lem: case1}
and Lemma \ref{Lem: case2}, we have the following result.
\begin{corollary}\label{Cor: case3}
Suppose $(S,V)$ is a marked surface with a decorated PH metric $(d,r)$, $\alpha\neq0$ is a constant
and $\overline{\mathbf{R}}: V\rightarrow(-\infty,2\pi)$ is a given function defined on the vertices.
\begin{description}
\item[(i)]
If $\alpha<0$ and there exists $h^\ast$ such that $\mathbf{R}_{\alpha,i}(h^\ast)<0$ for all $i\in V$,
then for any non-positive constant $\overline{\mathbf{R}}$,
there exists a decorated PH metric in $\mathcal{D}(d,r)$ with the combinatorial $\alpha$-curvature $\overline{\mathbf{R}}$;
\item[(ii)]
If $\alpha>0$ and there exists $h^\ast$ such that $\mathbf{R}_{\alpha,i}(h^\ast)>0$ for all $i\in V$,
then for any positive constant $\overline{\mathbf{R}}$,
there exists a decorated PH metric in $\mathcal{D}(d,r)$ with the combinatorial $\alpha$-curvature $\overline{\mathbf{R}}$.
\end{description}
\end{corollary}
\proof
\noindent\textbf{(i):}
Suppose there exists $h^\ast$ such that $\mathbf{R}_{\alpha,i}(h^\ast)<0$ for all $i\in V$.
By Lemma \ref{Lem: case1} (i),
for any constant combinatorial $\alpha$-curvature $\overline{\mathbf{R}}$ satisfying $\overline{\mathbf{R}}\leq\mathbf{R}_{\alpha,\min}(h^\ast)<0$,
there exists $\overline{h}$ with the constant combinatorial $\alpha$-curvature $\mathbf{R}_{\alpha,i}(\overline{h})
=\overline{\mathbf{R}}<0$ for all $i\in V$.
Taking $\overline{h}$ as the initial value of the combinatorial $\alpha$-Ricci flow with surgery (\ref{Eq: CRFS2}).
The conclusion follows from Lemma \ref{Lem: case1} (ii).

\noindent\textbf{(ii):}
The proof is similar to that of the case (i), so we omit it here.
\qed

The following result shows that the constant condition on $\overline{\mathbf{R}}$ in Corollary \ref{Cor: case3} could be removed.
\begin{theorem}\label{Thm: case4}
Suppose $(S,V)$ is a marked surface with a decorated PH metric $(d,r)$, $\alpha\neq0$ is a constant
and $\overline{\mathbf{R}}: V\rightarrow(-\infty,2\pi)$ is a given function defined on the vertices.
If there exists $h^\ast$ such that one of the following two conditions holds
\begin{description}
\item[(i)]
$\alpha<0,\ \overline{\mathbf{R}}_i\leq0,\ \mathbf{R}_{\alpha,i}(h^\ast)<0$ for all $i\in V$,
\item[(ii)]
$\alpha>0,\ \overline{\mathbf{R}}_i>0,\ \mathbf{R}_{\alpha,i}(h^\ast)>0$ for all $i\in V$,
\end{description}
then there exists a decorated PH metric in $\mathcal{D}(d,r)$ with the combinatorial $\alpha$-curvature $\overline{\mathbf{R}}$.
\end{theorem}
\proof
\noindent\textbf{(i):}
If $\overline{\mathbf{R}}$ is a constant, then the conclusion follows from Corollary \ref{Cor: case3} (i).
Suppose $\overline{\mathbf{R}}$ is not a constant,
set
\begin{equation*}
\overline{\mathbf{R}}_{\max}=\max_{i\in V}\overline{\mathbf{R}}_i,\
\overline{\mathbf{R}}_{\min}=\min_{i\in V}\overline{\mathbf{R}}_i.
\end{equation*}
By Corollary \ref{Cor: case3} (i), there exists $\overline{h}$ with the negative constant combinatorial $\alpha$-curvature $\mathbf{R}_{\alpha}(\overline{h})$ satisfying 
$\mathbf{R}_{\alpha}(\overline{h})<\overline{\mathbf{R}}_{\min}<0$.
Take $\overline{h}$ as the initial value of the combinatorial $\alpha$-Ricci flow with surgery (\ref{Eq: CRFS2}).
Then $M_{\alpha,i}(0)<0$ for all $i\in V$.
By Corollary \ref{Cor: application 1}, we have
$M_{\alpha,i}(t)\leq 0$ for all $i\in V$.
Submitting $\mathbf{A}_i<0$ and $M_{\alpha,i}(t)\leq0$ into (\ref{Eq: M alpha i}) gives
\begin{equation*}
\begin{aligned}
\frac{dM_{\alpha,i}}{dt}
&=\sum_{j\sim i}\frac{\mathbf{w}_{ij}}{e^{\alpha h_i}}(M_{\alpha,j}-M_{\alpha,i}) +\frac{\mathbf{A}_i}{e^{\alpha h_i}}M_{\alpha,i}-\alpha M_{\alpha,i}(M_{\alpha,i}+\overline{\mathbf{R}_i})\\
&\geq \sum_{j\sim i}\frac{\mathbf{w}_{ij}}{e^{\alpha h_i}}(M_{\alpha,j}-M_{\alpha,i})
-\alpha M_{\alpha,i}(M_{\alpha,i}+\overline{\mathbf{R}}_i)\\
&\geq \sum_{j\sim i}\frac{\mathbf{w}_{ij}}{e^{\alpha h_i}}(M_{\alpha,j}-M_{\alpha,i})
-\alpha M_{\alpha,i}(M_{\alpha,i}+\overline{\mathbf{R}}_{\max}).
\end{aligned}
\end{equation*}
By Theorem \ref{Thm: MX}, if $\mathbf{\overline{R}}_{\max}<0$, then
\begin{equation*}
0\geq M_{\alpha,i}(t) \geq
\frac{\mathbf{\overline{R}}_{\max}}
{-1+(\frac{\mathbf{\overline{R}}_{\max}}{M_{\alpha,\min}(0)}+1)
e^{\alpha \mathbf{\overline{R}}_{\max}t}}
\geq M_{\alpha,\min}(0)e^{-\alpha \mathbf{\overline{R}}_{\max}t}.
\end{equation*}
Then by (\ref{Eq: CRFS2}), we have
$M_{\alpha,\min}(0)e^{-\alpha\mathbf{\overline{R}}_{\max}t}\leq\frac{dh_i}{dt}
=M_{\alpha,i}(t)\leq0$.
This implies
\begin{equation}\label{Eq: F2}
\frac{M_{\alpha,\min}(0)}{\alpha \mathbf{\overline{R}}_{\max}}(1-e^{-\alpha \mathbf{\overline{R}}_{\max}t})\leq h_i(t)-h_i(0)\leq 0.
\end{equation}
By Theorem \ref{Thm: MX}, if $\mathbf{\overline{R}}_{\max}=0$, then
\begin{equation*}
0\geq M_{\alpha,i}(t)\geq \frac{M_{\alpha,\min}(0)}{\alpha M_{\alpha,\min}(0)t+1}.
\end{equation*}
This implies
$\frac{M_{\alpha,\min}(0)}{\alpha M_{\alpha,\min}(0)t+1}\leq\frac{dh_i}{dt}
=M_{\alpha,i}(t)\leq0$ by (\ref{Eq: CRFS2}).
Then
\begin{equation}\label{Eq: F3}
\frac{1}{\alpha}\ln (\alpha M_{\alpha,\min}(0)t+1)\leq h_i(t)-h_i(0)\leq0.
\end{equation}
The formulas (\ref{Eq: F2}) and (\ref{Eq: F3}) implies that $h_i(t)$ is uniformly bounded from above in $\mathbb{R}_{>0}$ for all $i\in V$.
By Lemma \ref{Lem: below 1},  $h_i(t)$ is uniformly bounded from below in $\mathbb{R}_{>0}$ for all $i\in V$ along the combinatorial $\alpha$-Ricci flow with surgery (\ref{Eq: CRFS2}).
By the arguments similar to that in Lemma \ref{Lem: case1} (ii), the solution of the combinatorial $\alpha$-Ricci flow with surgery (\ref{Eq: CRFS2}) with initial value $h(0)=\overline{h}$ exists for all time and converges.
The conclusion follows from Theorem \ref{Thm: equivalent 1} (i).

\noindent\textbf{(ii):}
If $\overline{\mathbf{R}}$ is a constant, then the conclusion follows from Corollary \ref{Cor: case3} (ii).
Suppose $\overline{\mathbf{R}}$ is not a constant.
By Corollary \ref{Cor: case3} (ii), there exists $\overline{h}$ with the positive constant combinatorial $\alpha$-curvature $\mathbf{R}_{\alpha}(\overline{h})$ satisfying $\mathbf{R}_{\alpha}(\overline{h})
>\overline{\mathbf{R}}_{\max}>0$.
Taking $\overline{h}$ as the initial value of the combinatorial $\alpha$-Ricci flow with surgery (\ref{Eq: CRFS2}).
Then $M_{\alpha,i}(0)>0$ for all $i\in V$.
By Corollary \ref{Cor: application 1}, we have
$M_{\alpha,i}(t)\geq 0$ for all $i\in V$.
Submitting $\mathbf{A}_i<0$ and $M_{\alpha,i}(t)\geq0$ into (\ref{Eq: M alpha i}) gives
\begin{equation*}
\begin{aligned}
\frac{dM_{\alpha,i}}{dt}
&=\sum_{j\sim i}\frac{\mathbf{w}_{ij}}{e^{\alpha h_i}}(M_{\alpha,j}-M_{\alpha,i}) +\frac{\mathbf{A}_i}{e^{\alpha h_i}}M_{\alpha,i}-\alpha M_{\alpha,i}(M_{\alpha,i}+\overline{\mathbf{R}}_i)\\
&\leq \sum_{j\sim i}\frac{\mathbf{w}_{ij}}{e^{\alpha h_i}}(M_{\alpha,j}-M_{\alpha,i})
-\alpha M_{\alpha,i}(M_{\alpha,i}+\overline{\mathbf{R}}_i)\\
&\leq \sum_{j\sim i}\frac{\mathbf{w}_{ij}}{e^{\alpha h_i}}(M_{\alpha,j}-M_{\alpha,i})
-\alpha M_{\alpha,i}(M_{\alpha,i}+\overline{\mathbf{R}}_{\min}).
\end{aligned}
\end{equation*}
Note that $\overline{\mathbf{R}}_{\min}>0$, by Theorem \ref{Thm: MX}, we have
\begin{equation*}
0\leq M_{\alpha,i}(t)
\leq \frac{\mathbf{\overline{R}}_{\min}}
{-1+(\frac{\mathbf{\overline{R}}_{\min}}{M_{\alpha,\max}(0)}+1)
e^{\alpha \mathbf{\overline{R}}_{\min}t}}
\leq M_{\alpha,\max}(0)e^{-\alpha \mathbf{\overline{R}}_{\min}t}.
\end{equation*}
Then by (\ref{Eq: CRFS2}), we have
$0\leq \frac{dh_i}{dt}=M_{\alpha,i}(t)\leq M_{\alpha,\max}(0)e^{-\alpha \mathbf{\overline{R}}_{\min}t}$.
Hence,
\begin{equation*}
0\leq h_i(t)-h_i(0)\leq \frac{M_{\alpha,\max}(0)}{\alpha \mathbf{\overline{R}}_{\min}}(1-e^{-\alpha \mathbf{\overline{R}}_{\min}t}).
\end{equation*}
This implies that the solution of the combinatorial $\alpha$-Ricci flow with surgery (\ref{Eq: CRFS2}) with initial value $h(0)=\overline{h}$ lies in a compact subset of $\mathbb{R}_{>0}^N$.
Combining with $\frac{dh_i}{dt}\geq 0$,
the solution of the combinatorial $\alpha$-Ricci flow with surgery (\ref{Eq: CRFS2}) with initial value $h(0)=\overline{h}$ exists for all time and converges.
The conclusion follows from Theorem \ref{Thm: equivalent 1} (i).
\qed

\noindent\textbf{Proof of Theorem \ref{Thm: main 1}.}
\noindent\textbf{(i):}
Since $\chi(S)<0$, by Theorem \ref{Thm: BL},
there exists $\overline{h}$ such that $\sum^N_{i=1} K_i(\overline{h})>2\pi \chi(S)$ and $K_i(\overline{h})<0$ for all $i\in V$.
Then $\mathbf{R}_{\alpha,i}(\overline{h})
=\frac{\mathbf{K}_i(\overline{h})}{e^{\alpha \overline{h}_i}}<0$ for all $i\in V$.
By Theorem \ref{Thm: case4} (i), there exists a decorated PH metric in $\mathcal{D}(d,r)$ with the non-positive combinatorial $\alpha$-curvature $\overline{\mathbf{R}}$.

\noindent\textbf{(ii):}
By Theorem \ref{Thm: BL},
there exists $\overline{h}$ such that $\sum^N_{i=1} K_i(\overline{h})>2\pi \chi(S)$ and $K_i(\overline{h})>0$ for all $i\in V$.
Then $\mathbf{R}_{\alpha,i}(\overline{h})
=\frac{\mathbf{K}_i(\overline{h})}{e^{\alpha \overline{h}_i}}>0$ for all $i\in V$.
By Theorem \ref{Thm: case4} (ii), there exists a decorated PH metric in $\mathcal{D}(d,r)$ with the positive combinatorial $\alpha$-curvature $\overline{\mathbf{R}}$.

\noindent\textbf{(iii):}
The conclusion follows from Theorem \ref{Thm: BL}.
\qed

Combining Theorem \ref{Thm: equivalent 1} (ii) and Theorem \ref{Thm: main 1},
we have the following longtime existence and convergence for the solution of the combinatorial $\alpha$-Ricci flow with surgery (\ref{Eq: CRFS}),
which corresponds to the case of the combinatorial $\alpha$-Ricci flow with surgery (\ref{Eq: CRFS}) in Theorem \ref{Thm: main 2}.

\begin{corollary}
Suppose $(S,V)$ is a marked surface with a decorated PH metric $(d,r)$, $\alpha\in \mathbb{R}$ is a constant and $\overline{\mathbf{R}}: V\rightarrow (-\infty,2\pi)$ is a given function defined on the vertices.
If one of the conditions (i)(ii)(iii) in Theorem \ref{Thm: main 1} is satisfied,
then the solution of the combinatorial $\alpha$-Ricci flow with surgery (\ref{Eq: CRFS}) exists for all time and converges exponentially fast to $\overline{h}$ such that the decorated PH metric $(d(\overline{h}),r(\overline{h}))$ has the combinatorial $\alpha$-curvature $\overline{\mathbf{R}}$.
\end{corollary}

\section{Combinatorial $\alpha$-Calabi flow with surgery}
\label{Sec: CCFS}

Analogy to Definition \ref{Def: CRFS} for the combinatorial $\alpha$-Ricci flow with surgery,
we introduce the following combinatorial $\alpha$-Calabi flow with surgery.

\begin{definition}\label{Def: CCFS}
Suppose $(S,V)$ is a marked surface with a decorated PH metric $(d,r)$, $\alpha\in \mathbb{R}$ is a constant and $\overline{\mathbf{R}}: V\rightarrow\mathbb{R}$ is a given function defined on the vertices.
The combinatorial $\alpha$-Calabi flow with surgery is defined to be
\begin{eqnarray}\label{Eq: CCFS}
\begin{cases}
\frac{dh_i}{dt}
=\Delta_\alpha(\overline{\mathbf{R}}
-\mathbf{R}_\alpha)_i,\\
h_i(0)=h_0.
\end{cases}
\end{eqnarray}
\end{definition}

Parallelling to Theorem \ref{Thm: equivalent 1}, we have the following result on the combinatorial $\alpha$-Calabi flow with surgery (\ref{Eq: CCFS}).

\begin{theorem}\label{Thm: equivalent 2}
Suppose $(S,V)$ is a marked surface with a decorated PH metric $(d,r)$, $\alpha\in \mathbb{R}$ is a constant and $\overline{\mathbf{R}}: V\rightarrow \mathbb{R}$ is a given function defined on the vertices.
\begin{description}
\item[(i)]
If the solution of the combinatorial $\alpha$-Calabi flow with surgery (\ref{Eq: CCFS}) converges,
then there exists a decorated PH metric in $\mathcal{D}(d,r)$ with the combinatorial $\alpha$-curvature $\overline{\mathbf{R}}$.
\item[(ii)]
If $\overline{\mathbf{R}}: V\rightarrow(-\infty,2\pi)$ satisfies $\alpha\overline{\mathbf{R}}\geq0$ and
there exists a decorated PH metric in $\mathcal{D}(d,r)$ with the combinatorial $\alpha$-curvature $\overline{\mathbf{R}}$,
then the solution of the combinatorial $\alpha$-Calabi flow with surgery (\ref{Eq: CCFS}) exists for all time and converges exponentially fast for any initial $h(0)\in \mathbb{R}_{>0}^N$.
\end{description}
\end{theorem}
\proof
\noindent\textbf{(i):}
The proof is almost the same as that of Theorem \ref{Thm: local convergence} (i), so we omit it here.
%Suppose the solution $h(t)$ of the combinatorial $\alpha$-Calabi flow with surgery (\ref{Eq: CCFS}) converges to $\overline{h}$ as $t\rightarrow +\infty$, then by the $C^1$-smoothness of $\mathbf{R}_\alpha$, we have $\mathbf{R}_\alpha(\overline{h})=\lim_{t\rightarrow +\infty}\mathbf{R}_\alpha(h(t))$.
%Furthermore, there exists a sequence $t_n\in(n,n+1)$ such that
%\begin{equation*}
%h_i(n+1)-h_i(n)=h'_i(t_n)
%=\Delta_\alpha(\mathbf{\overline{R}}-
%\mathbf{R}_{\alpha}(h(\xi_n)))_i\rightarrow 0,\  \text{as}\  n\rightarrow +\infty.
%\end{equation*}
%By Lemma \ref{Lem: matrix negative} and Definition \ref{Def: Delta alpha},
%the Laplace operator $\Delta_\alpha$ is negative definite.
%Then $\mathbf{R}_{\alpha,i}(\overline{h})
%=\lim_{n\rightarrow +\infty}\mathbf{R}_{\alpha,i}(h(t_n))
%=\overline{\mathbf{R}}_i$ for all $i\in V$
%and $\overline{h}$ is a discrete conformal factor with the combinatorial $\alpha$-curvature $\overline{\mathbf{R}}$.

\noindent\textbf{(ii):}
By the arguments similar to that in Theorem \ref{Thm: equivalent 1} (ii),
the function $\mathcal{W}(h)$ defined by (\ref{Eq: W alpha}) can also be applied to this theorem.
Along the combinatorial $\alpha$-Calabi flow with surgery (\ref{Eq: CCFS}), we have
\begin{align*}
\frac{d\mathcal{W}_\alpha(h(t))}{dt}
&=\sum^N_{i=1}\frac{\partial \mathcal{W}_\alpha}{\partial h_i}\frac{dh_i}{dt}
=\sum^N_{i=1}(\mathbf{R}_{\alpha}
-\overline{\mathbf{R}})_ie^{\alpha h_i}
\Delta_\alpha(\overline{\mathbf{R}}
-\mathbf{R}_{\alpha})_i\\
&=-(\mathbf{R}_\alpha-\overline{\mathbf{R}})^\mathrm{T}\cdot \mathbf{L}\cdot (\mathbf{R}_\alpha-\overline{\mathbf{R}})\geq 0
\end{align*}
by (\ref{Eq: F67}), where the last inequality follows from the negative definiteness of $\mathbf{L}$ in Lemma \ref{Lem: matrix negative}.
This implies $\mathcal{W}_\alpha(h(0))\leq \mathcal{W}_\alpha(h(t))\leq 0$.
Combining this with the properness of $\mathcal{W}_\alpha(h)$,
the solution $h(t)$ of the combinatorial $\alpha$-Calabi flow with surgery (\ref{Eq: CCFS}) is uniformly bounded from above in $\mathbb{R}^N_{>0}$.
The following Lemma \ref{Lem: below 2} shows that $h(t)$ is uniformly bounded from below in $\mathbb{R}^N_{>0}$.
Hence, the solution $h(t)$ of the combinatorial $\alpha$-Calabi flow with surgery (\ref{Eq: CCFS}) lies in a compact subset of $\mathbb{R}^N_{>0}$.
This implies the solution of the combinatorial $\alpha$-Calabi flow with surgery (\ref{Eq: CCFS}) exists for all time and $\lim_{t\rightarrow +\infty}\mathcal{W}_\alpha(h(t))$ exists.
Moreover, there exists a sequence $t_n\in(n,n+1)$ such that as $n\rightarrow +\infty$,
\begin{equation*}
\begin{aligned}
&\mathcal{W}_\alpha(h(n+1))-W_\alpha(h(n))
=(\mathcal{W}_\alpha(h(t))'|_{t_n}
=\nabla \mathcal{W}_\alpha\cdot\frac{dh_i}{dt}|_{t_n}\\
=&\sum^N_{i=1}(\mathbf{R}_{\alpha}
-\overline{\mathbf{R}})_ie^{\alpha h_i}
\Delta_\alpha(\overline{\mathbf{R}}
-\mathbf{R}_{\alpha})_i|_{t_n}
=-(\mathbf{R}_\alpha-\overline{\mathbf{R}})^\mathrm{T}
\cdot \mathbf{L}\cdot (\mathbf{R}_\alpha-\overline{\mathbf{R}})|_{t_n}
\rightarrow 0.
\end{aligned}
\end{equation*}
Then $\lim_{n\rightarrow +\infty}\mathbf{R}_{\alpha,i}(h(t_n))
=\overline{\mathbf{R}}_i
=\mathbf{R}_{\alpha,i}(\overline{h})$ for all $i\in V$.
By $\{h(t)\}\subset\subset \mathbb{R}^N_{>0}$,
there exists $h^*\in \mathbb{R}^N_{>0}$ and a convergent subsequence of $\{h(t_n)\}$, still denoted as $\{h(t_n)\}$ for simplicity, such that $\lim_{n\rightarrow \infty}h(t_n)=h^*$.
This implies
$\mathbf{R}_\alpha(h^*)=\lim_{n\rightarrow +\infty}\mathbf{R}_\alpha(h(t_n))
=\mathbf{R}_\alpha(\overline{h})$.
Then $h^*=\overline{h}$ by Theorem \ref{Thm: rigidity}.
Therefore, $\lim_{n\rightarrow \infty}h(t_n)=\overline{h}$.

Similar to the proof of Theorem \ref{Thm: local convergence} (ii),
set $\Gamma(h)=\Delta_\alpha(\overline{\mathbf{R}}
-\mathbf{R}_\alpha)$,
then $D\Gamma|_{h=\overline{h}}$ has $N$ negative eigenvalues.
This implies that $\overline{h}$ is a local attractor of the combinatorial $\alpha$-Calabi flow with surgery (\ref{Eq: CCFS}).
The conclusion follows from Lyapunov Stability Theorem (\cite{Pontryagin}, Chapter 5).
\qed

\begin{lemma}\label{Lem: below 2}
Under the same assumptions in Lemma \ref{Thm: equivalent 2} (ii),
along the combinatorial $\alpha$-Calabi flow with surgery (\ref{Eq: CCFS}),
$h_i(t)$ is uniformly bounded from below in $\mathbb{R}_{>0}$ for all $i\in V$.
\end{lemma}
\proof
Similar to the arguments in the proof of Lemma \ref{Lem: below 1},
we just need to prove this lemma for some fixed weighted Delaunay triangulation $\mathcal{T}$ in $(d,r)$.
Direct calculations give
\begin{equation}\label{Eq: F4}
\begin{aligned}
e^{\alpha h_i}\cdot\frac{d h_i}{dt}
&=e^{\alpha h_i}\cdot\Delta^{\mathcal{T}}_\alpha
(\overline{R}-R_\alpha)_i\\
&=\frac{\partial K_i}{\partial h_i}(\overline{R}_i-R_{\alpha,i})
+\sum_{j\neq i,\ j\in V}\frac{\partial K_i}{\partial h_j}(\overline{R}_j-R_{\alpha,j})\\
&=-\sum_{j\sim i}w_{ij}\cosh l_{ij}(\overline{R}_i- R_{\alpha,i})
+\sum_{j\sim i}w_{ij}(\overline{R}_j-R_{\alpha,j})\\
&=\sum_{j\sim i}w_{ij}[\cosh l_{ij}(R_{\alpha,i}-\overline{R}_i)+(\overline{R}_j-R_{\alpha,j})],
\end{aligned}
\end{equation}
where the third line use (\ref{Eq: A i}) and (\ref{Eq: L B}).
Note that in the proof of Theorem \ref{Thm: equivalent 2} (ii),
the solution $h(t)$ of the combinatorial $\alpha$-Calabi flow with surgery (\ref{Eq: CCFS}) is uniformly bounded from above in $\mathbb{R}^N_{>0}$,
denoted by $C_7$, i.e., $h_i(t)\leq C_7$ for all $i\in V$.
Then $e^{\alpha h_i}>1$ for $\alpha\geq0$ and $e^{\alpha h_i}\geq e^{\alpha C_7}$ for $\alpha<0$.
By the definition of $K_i$, we have $(2-N)\pi<K_i<2\pi$ for all $i\in V$.
Then $|R_{\alpha,i}|=|\frac{K_i}{e^{\alpha h_i}}|<(2+N)\pi$ for $\alpha\geq0$ and $|R_{\alpha,i}|<\frac{(2+N)\pi}{e^{\alpha C_7}}$ for $\alpha<0$.
Hence for $\alpha\in \mathbb{R}$, $R_{\alpha,i}$ is bounded for all $i\in V$.

By contradiction, we assume that there exists at least one vertex $i\in V$ such that $\lim_{t\rightarrow T}h_i(t)=0$ for $T\in(0,+\infty]$.
By the proof of Lemma \ref{Lem: below 1},
there exists a positive constant $C_5$ such that if $h_i<C_5$, then $R_{\alpha,i}-\overline{R}_i>0$.
By the map (\ref{Eq: h,r}), we have $\lim_{t\rightarrow T}r_i(t)=+\infty$.
Combining with (\ref{Eq: inversive distance}), we have
$\cosh l_{ij}\rightarrow +\infty$.
Hence, we can choose a sufficient small $C_8<C_5$ such that if $h_i<C_8$, then $R_{\alpha,i}-\overline{R}_i>0$ and
\begin{equation}\label{Eq: F7}
\cosh l_{ij}(R_{\alpha,i}-\overline{R}_i)
+(\overline{R}_j-R_{\alpha,j})
>\cosh l_{ij}(R_{\alpha,i}-\overline{R}_i)
-|\overline{R}_j|-|R_{\alpha,j}|
>0.
\end{equation}
Note that $\sum_{j\sim i}w_{ij}>0$ by Lemma \ref{Lem: L decomposition}, submitting (\ref{Eq: F7}) into (\ref{Eq: F4}) gives $e^{\alpha h_i}\cdot\frac{d h_i}{dt}>0$.
Hence $\frac{d h_i}{dt}
=\Delta^{\mathcal{T}}_\alpha
(\overline{R}-R_\alpha)_i>0$.
The rest of the proof is similar to that of Lemma \ref{Lem: below 1}, so we omit it here.
\qed

Combining Theorem \ref{Thm: equivalent 2} (ii) and Theorem \ref{Thm: main 1},
we have the following longtime existence and convergence for the solution of the combinatorial $\alpha$-Calabi flow with surgery (\ref{Eq: CCFS}),
which corresponds to the case of the combinatorial $\alpha$-Calabi flow with surgery (\ref{Eq: CCFS}) in Theorem \ref{Thm: main 2}.

\begin{corollary}
Suppose $(S,V)$ is a marked surface with a decorated PH metric $(d,r)$, $\alpha\in \mathbb{R}$ is a constant and $\overline{\mathbf{R}}: V\rightarrow (-\infty,2\pi)$ is a given function defined on the vertices.
If one of the conditions (i)(ii)(iii) in Theorem \ref{Thm: main 1} is satisfied,
then the solution of the combinatorial $\alpha$-Calabi flow with surgery (\ref{Eq: CCFS}) exists for all time and converges exponentially fast to $\overline{h}$ such that the decorated PH metric $(d(\overline{h}),r(\overline{h}))$ has the combinatorial $\alpha$-curvature $\overline{\mathbf{R}}$.
\end{corollary}

\end{document}